\documentclass{article}
\usepackage[utf8]{inputenc}
\title{\normalsize{12th AIMMS-MOPTA Optimization Modeling Competition}\\
\Large{Vehicle Fleet Sizing, Positioning and Routing Problem with Stochastic Customers}}
\author{Team M$\alpha$D Condors\\
Maria J. Consuegra-Laino, Alfaima L. Solano-Blanco, David Corredor-Montenegro\\Advisor: Camilo Gómez}
\date{%
    \textit{Departamento de Ingeniería Industrial}, Universidad de los Andes\\[2ex]%
    }

\usepackage{natbib}
\usepackage{amsmath}
\usepackage{enumerate}
\usepackage{dsfont}
\usepackage{comment}
\usepackage{geometry}
\DeclareMathOperator*{\argmin}{arg\,min}

\usepackage{bm}
\usepackage{amssymb}
\usepackage{amsfonts}
\usepackage{amsthm}
\usepackage{manfnt}
\usepackage{mathtools}
\usepackage{soul}
\usepackage{mathtools}
\usepackage{tikz}
\usepackage{pgfplots}
\usepgfplotslibrary{fillbetween}
\usepackage{verbatim}
\usepackage{mathrsfs}
\usepackage{algorithmicx}
\usepackage{algpseudocode}
\usepackage{algorithm}
\usepackage{natbib}
\usepackage{subcaption}
\usepackage{verbatim}

\usepackage{caption}
\usepackage{graphicx}
\usepackage{booktabs}

\usetikzlibrary{shapes,arrows,positioning}
\usetikzlibrary{calc, automata, chains, arrows.meta}
\usetikzlibrary{math}
\usetikzlibrary{matrix, positioning, fit}
\numberwithin{equation}{section}

\usepackage{blindtext}

\DeclarePairedDelimiter{\ceil}{\lceil}{\rceil}

\usepackage{graphicx}
\graphicspath{ {./images/} }

\newcommand{\Mtodo}{\textcolor{red}}   
\newcommand{\Atodo}{\textcolor{violet}}   
\newcommand{\Dtodo}{\textcolor{blue}}   
\newcommand{\Ctodo}{\textcolor{orange}}   

\usepackage{optidef}

\newcommand{\G}{\mathcal{G}}
\newcommand{\N}{\mathcal{N}}
\newcommand{\A}{\mathcal{A}}
\newcommand{\D}{\mathcal{D}}
\newcommand{\Rout}{\mathscr{R}}

\newcommand{\Scen}{\mathcal{S}}
\newcommand{\Hh}{\mathcal{H}}

\definecolor{ao(english)}{rgb}{0.0, 0.5, 0.0}

\providecommand{\keywords}[1]{\textbf{\textit{Keywords: }} #1}

\begin{document}
\maketitle

\begin{abstract}
The Vehicle Fleet Sizing, Positioning and Routing Problem with Stochastic Customers (VFSPRP-SC) consists on pairing strategic decisions of depot positioning and fleet sizing with operational vehicle routing decisions while taking into account the inherent uncertainty of demand. We successfully solve the VFSPRP-SC with a methodology comprised of two main blocks: i) a scenario generation phase and ii) a two--stage stochastic program. For the first block, a set of scenarios is selected with a simulation--based approach that captures the behavior of the demand and allows us to come up with different solutions that could match different risk profiles. The second block is comprised of a facility location and allocation model and a Multi Depot Vehicle Routing Problem (MDVRP) assembled under a two--stage stochastic program. We propose several novel ideas within our methodology: problem specific cuts that serve as an approximation of the expected second stage costs as a function of first stage decisions; an activation paradigm that guides our main optimization procedure; and, a way of mapping feasible routes from one second--stage problem data into another; among others. We performed experiments for two cases: the first case considers the expected value of the demand, and the second case considers the \textit{right tail} of the demand distribution, seeking a conservative solution. For the evaluation of the solutions, we carried out an ex-post analysis by running the routing model for the demand days available in the 12th AIMMS-MOPTA Optimization Modeling Competition dataset. The conservative solution locates 5 depots and allocates 21 vehicles. With the ex-post analysis, we provide relevant metrics like costs, vehicle utilization, service level (of 100\% in the conservative case) and average daily emissions, among others, for the decision-makers to compare the solutions as alternatives. By using acceleration techniques throughout the methodology we obtain solutions within 1 to 6 hours, reasonable times considering the strategic nature of the decision. For the ex-post evaluation, we run 75\% of the instances in less than 3 minutes, so the methodology used to solve the MDVRP is well suited for daily operation. We implemented our solution scheme in AIMMS as a user-friendly decision support system.

\end{abstract}

\keywords{Vehicle Routing Problem, Multi-Depot; Fleet Sizing; Facility location; Stochastic Customers; Time windows; Scenario generation; L-shaped method; Pulse algorithm}

\section{Introduction}\label{s:ProblemDefinition}

Optimal distribution in supply chains is relevant as it allows savings for companies and consumers. A crucial challenge in supply chain design and operation arises from the uncertainty in external factors, such as demand. Businesses are in the growing need for robust systems that can react and adapt to the uncertain and ever-changing business environment. Thus, the layout of the distribution structure must guarantee the cutback of logistics costs and an adequate service level. The goal of this research paper is to solve an optimal vehicle routing problem where in addition, the location of the depots for the vehicles and the fleet size should be optimized to comply with uncertain demand at minimum cost.

Mathematical models and effective solution techniques to face location and routing problems can be found in literature. The Vehicle Routing Problem (VRP) consists on determining an optimal group of routes to be performed by a fleet of vehicles to serve a given group of customers. The VRP is one of the most important combinatorial optimization problems, and has been established as an NP-hard problem along with its variants. Some comprehensive surveys about the variants and solution approaches are \cite{cordeau2007vehicle,toth2014vehicle} for exact models and \cite{laporte2002classical, cordeau2002guide} for heuristics and metaheuristics.

One common approach for tackling the VRP is to create a pool of routes that comply with the necessary constraints with an auxiliary model, which are subsequently selected by an assignment model. Traditionally the problem solved by the auxiliary model is reduced to constrained shortest path problems (CSP). Many approaches have been proposed in literature to tackle the CSP, including dynamic programming based labeling algorithms \citep{Dumitrescu2003,Joksch1966,Thomas2019} and path ranking approaches \citep{Handler1980,Santos2007,Sedeno-Noda2015}. A competitive approach for the CSP was proposed by \cite{Lozano2013} and coined with the term \textit{pulse algorithm}. This exact algorithm uses depth-first search combined with effective pruning strategies to make an implicit exploration of the solution space of the CSP. In addition, the pulse algorithm has been successfully extended to solve other hard shortest path variants, such as the elementary shortest path problem with resource constraints \citep{Lozano2016,Li2019}, and the weight constrained shortest path problem with replenishment \citep{Bolivar2014}, among other shortest path variants.

In some contexts, it is important to consider variability in the problem parameters. Some variants of the VRP that consider uncertainty are: the Stochastic demand VRP, where the demand is only known when the vehicle arrives to its destination; the VRP with stochastic times, where the traversing times through the network are treated as random variables; and, the VRP with stochastic customers, in which each node has a certain probability of placing an order. For the stochastic customers VRP, the common solution schemes consist in generating valid routes in a first stage where the realizations are unknown and using those routes in a second stage given the realizations of present and absent customers \citep{gendreau2014chapter}. The problem has been solved with semifixed routes \citep{waters1989vehicle, benton1992vehicle, sungur2010model}, L-shaped methods for up to 46 nodes \citep{gendreau1995exact}, heuristics as tabu search \citep{gendreau1996tabu} and metaheuristics \citep{balaprakash2015estimation}. 

At the strategic level for a delivery company, decisions on the locations that are most suitable for the execution of the daily operation are also affected by stochasticity \citep{snyder2006facility}.  A variant of the VRP that includes the location decisions and the fleet size decisions have been also addressed in the literature \citep{kocc2016fleet,golden1984fleet,wu2002heuristic}. However, these authors did not consider uncertainty which is one of the main challenges of the research.

In this research, we address a problem with the combination of several of the mentioned variants of the VRP: the Vehicle Fleet Sizing, Positioning and Routing Problem with Stochastic Customers (henceforth VFSPRP-SC). This problem consists on finding a subset of locations to install depots for which a number of vehicles must be determined to guarantee demand satisfaction. The vehicles assigned to each depot will remain in place for the entire duration of the problem, even if they have to remain idle on certain days. Finally, the vehicles need to be routed over a transportation network to satisfy the stochastic daily demand of the active clients.
\begin{figure}
    \centering
    \includegraphics[width=0.7\textwidth]{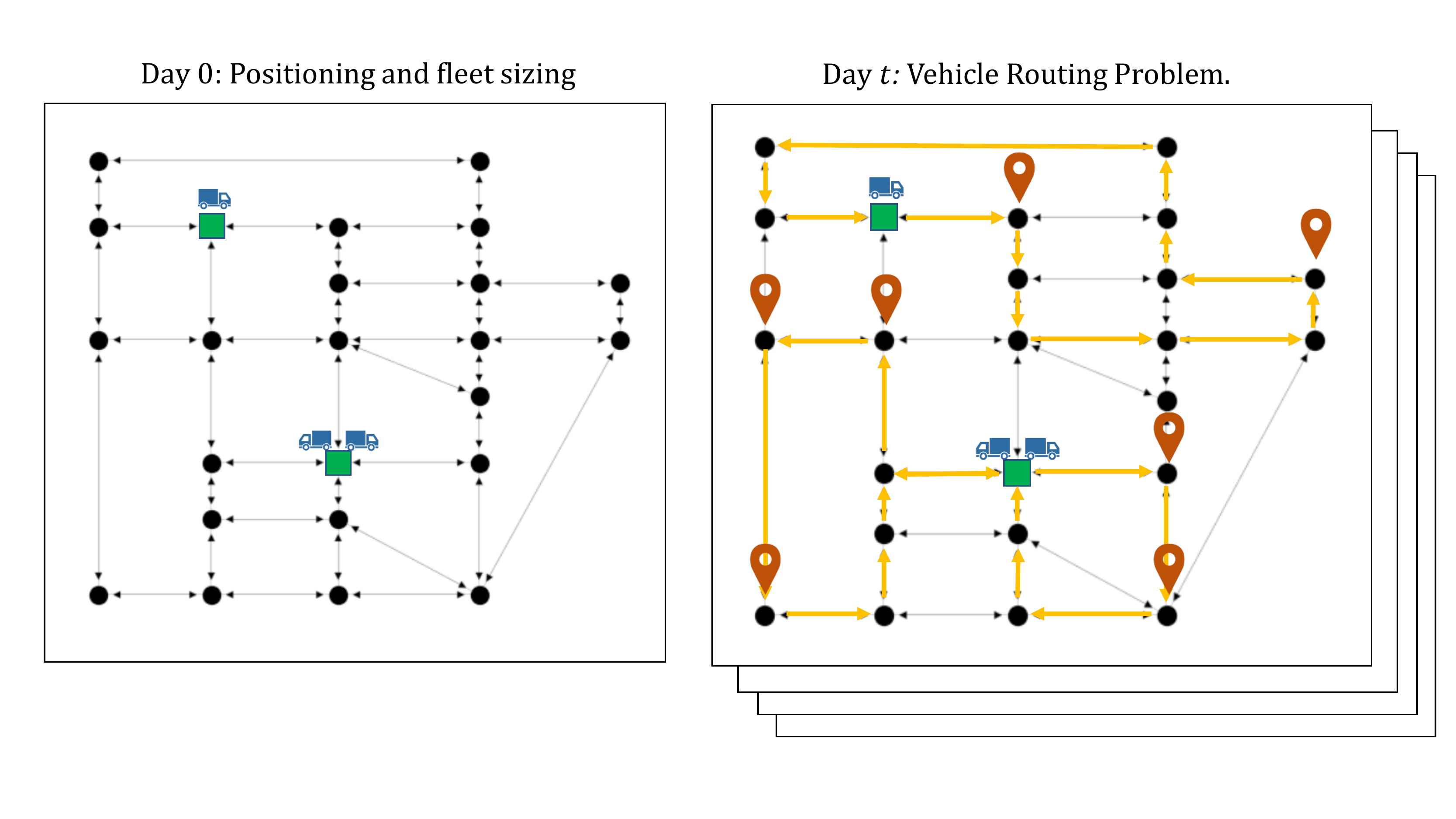}
    \caption{Visual representation of the Vehicle Fleet Sizing, Positioning and Routing Problem with Stochastic Customers (VFSPRP-SC) where the black dots represent locations, the green squares represent installed depots, the red marks represent the active clients and the yellow lines represent the route of vehicles.}
    \label{fig:probIlustration}
\end{figure}

\begin{figure}[h]
    \centering
    \begin{tikzpicture}[font=\ttfamily,
        mymatrix/.style={matrix of nodes, nodes=typetag, row sep=1em},
        mycontainer/.style={draw=black,dashed,very thick, inner sep=1ex},
        mycontainerFS/.style={draw=ao(english),dashed,thick, inner sep=1ex},
        mycontainerSS/.style={draw=blue,dashed,thick, inner sep=1ex},
        mycontainerMDVRP/.style={draw=orange,dashed,thick, inner sep=1ex},
        mycontainerSimul/.style={draw=red,dashed, thick, inner sep=1ex},
        typetag/.style={circle,draw=gray, inner sep=1ex, anchor=west},
        title/.style={draw=none, color=black, inner sep=0pt},
        block/.style={rectangle, draw, fill=blue!20, 
        text width=5em, text centered, rounded corners, minimum height=3em},
        connect/.style={circle, draw, fill=yellow!20, minimum height=1em},
        line/.style={draw, -latex'},node distance = 2cm, auto]

      \node[block](FL) {Facility location};
      \node[connect,below =1.5cm of FL](c1) {\small $1$};
      \node[block, below left=1cm and 2.5cm of c1 ](MDVRP1) {Set covering$_1$};
      \node[block, below right =1cm and 2.5cm of c1 ](MDVRP2) {Set covering$_n$};
      \node[connect,below =1.5cm of MDVRP1](c2) {\small $2$};
      \node[connect,below =1.5cm of MDVRP2](c3) {\small $2$};

      \node[block, below left =1cm and 0.7cm of c2 ](pul1) {$\text{ESPPRC-R}_{1,1}$};
      \node[block, below right =1cm and 0.7cm of c2 ](pul2) {$\text{ESPPRC-R}_{1,h}$};
      
      \node[block, below left =1cm and 0.7cm of c3 ](pul3) {$\text{ESPPRC-R}_{n,1}$};
      \node[block, below right =1cm and 0.7cm of c3 ](pul4) {$\text{ESPPRC-R}_{n,h}$};
      
      \node[block, below = 8cm of c1] (simul) {Simulation};
      \node[block, left = 2.5cm of simul] (Data) {Data analysis};
      \node[block, right = 2.5cm of simul] (select) {Scenario selection};

      \node[above left = 0.5 cm and 0 cm of c1](x){Solution};
      \node[above right = 0.5 cm and 0.3 cm of c1](cut){Cuts};
      
      \node[above left = 0.7 cm and 0 cm of MDVRP2]{Feedback};
      \node[above right = 0.7 cm and 0 cm of MDVRP1]{Feedback};
      
      \node[above left = 0.5 cm and 0 cm of c3](dual){Dual info};
      \node[above right = 0.5 cm and 0.3 cm of c3](rout){Routes};
      
      \node[above left = 0.5 cm and 0 cm of c2](dual1){Dual info};
      \node[above right = 0.5 cm and 0.3 cm of c2](rout1){Routes};
      
      \node[below= 3.5cm of c1](dots){\large $\bm{\cdots}$};
      \node[below= 1.2cm of c2]{\large $\bm{\cdots}$};
      \node[below= 1.2cm of c3]{\large $\bm{\cdots}$};
      

      \node[below left =0.6 cm and 1cm of pul1](a1) {};
      \node[below right =0.6cm and 1cm of pul4](a2) {};
      
      \node[below left =0.4cm and 0.5cm of pul1](ass1) {};
      \node[below right =0.4cm and 0.5 cm of pul4](ass2) {};
      \node[below left = 0.6cm and 6cm of c1](ass3) {};
      \node[below right = 0.6cm and 6cm of c1](ass4) {};

      \node[below left =0.2cm and 0cm of pul1](a3) {};
      
      \node[below right =0.2cm and 0cm of pul4](a4) {};
      
      \node[above= 0.5cm of FL](a0) {};
      \node[above left = 0cm and 6.1 cm of FL](afs1) {};
      \node[above right = 0cm and 6.1cm of FL](afs2) {};
      \node[above left = 0.4cm and 6.1cm of c1](afs3) {};
      \node[above right = 0.4cm and 6.1cm of c1](afs4) {};
      
      \node[below right =0.6 and 0.2 of a1](aScenL1){};
      \node[below of= aScenL1](aScenL2){};
      \node[below left =0.6 and 0.2 of a2](aScenR1){};
      \node[below of= aScenR1](aScenR2){};

      \node[below left = 0.2 and 0 of aScenL2](fin1){};
      \node[below right = 0.2 and 0 of aScenR2](fin2){};
      
        \node[mycontainer,fit= (a0) (fin1) (fin2)](cont1){};
       \node[anchor=south west] at (cont1.north west) {VFSPRP-SC (See \S \ref{s:SolMethodology})};
       
       \node[mycontainerMDVRP,fit= (MDVRP1) (a3) (pul2)](cont2){};
       \node[anchor=south west] at (cont2.north west) {\footnotesize MDVRP};
       
       \node[mycontainerMDVRP,fit= (pul3) (MDVRP2) (a4) ](cont6){};

       \node[mycontainerFS,fit= (afs1) (afs2) (afs3) (afs4)](cont3){};
       \node[anchor=south east] at (cont3.north east) {\footnotesize First Stage (See \S \ref{ss:LShaped})};
       
       \node[mycontainerSS,fit= (ass1) (ass2) (ass3) (ass4)](cont4){};
       \node[anchor=south east] at (cont4.north east) {\footnotesize Second Stage (See \S \ref{ss:SecondStage})};
       
       \node[mycontainerSimul,fit= (aScenL1) (aScenL2) (aScenR1) (aScenR2)](cont5){};
       \node[anchor=south east] at (cont5.north east) {\footnotesize Scenario generation (See \S \ref{ss:Simulations})};

      \path [line] (FL) -- (x)-- (c1);
      \path [line] (c1) -- (cut)-- (FL);
      
      \path [line] ([xshift=-1cm] c1) |- (MDVRP1); 
      \path [line] ([xshift=1cm] c1) |- (MDVRP2); 

      \path [line] (MDVRP1) |- (c1);
      \path [line] (MDVRP2) |- (c1);
      
     
      \path [line] (MDVRP1) -- (dual1)-- (c2);
      \path [line] (c2) -- (rout1)-- (MDVRP1);
      
      \path [line] (MDVRP2) -- (dual)-- (c3);
      \path [line] (c3) -- (rout)-- (MDVRP2);

      \path [line] ([xshift=-1cm] c2) |- ([yshift=.5cm] pul1);
      \path [line] ([xshift=1cm] c2) |- ([yshift=.5cm] pul2);
      
      \path [line] ([xshift=-1cm]c3) |- ([yshift=.5cm] pul3);
      \path [line] ([xshift=1cm]c3) |- ([yshift=.5cm] pul4);
      
      \path [line] (pul1) |- (c2) ;
      \path [line] (pul2) |- (c2) ;
      
      \path [line] (pul3) |- (c3) ;
      \path [line] (pul4) |- (c3) ;
      
      \path[line] (Data)|-(simul);
      \path[line](simul)|-(select);

      \path[line]([xshift=-.7cm]cont5)|-([yshift=-0.4cm]cont2);
      \path[line]([xshift=.7cm]cont5)|-([yshift=-0.4cm]cont6);
      
    \end{tikzpicture}
    \caption{Solution methodology for the VFSPRP-SC with the first stage  (outlined in green), the second stage (outlined in blue), each MDVRP subproblem (outlined in orange), and the scenario generation scheme (outlined in red).}
    \label{fig:Methodology}
\end{figure}
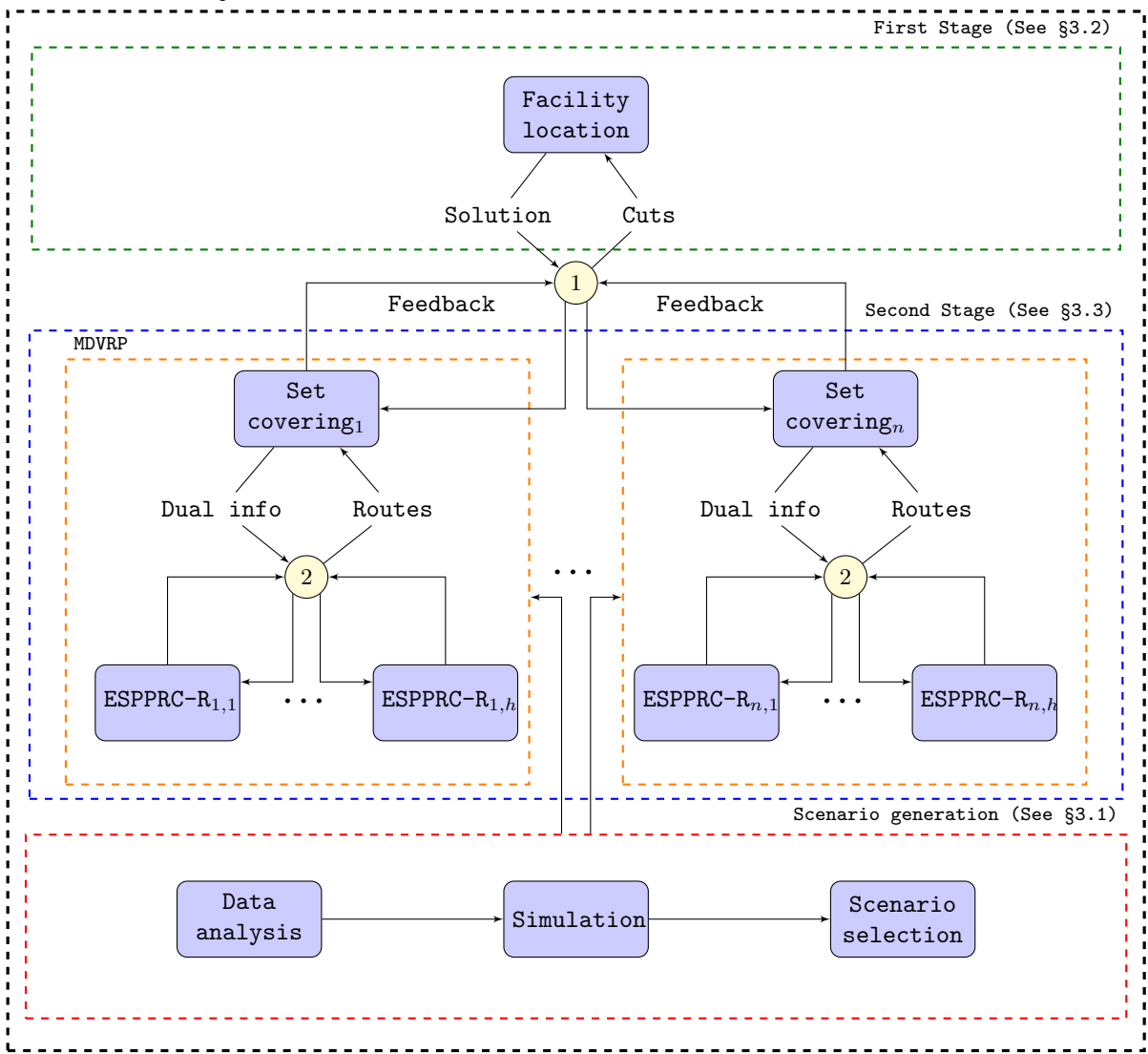

We propose a solution approach for the VFSPRP-SC that uses various interconnected models as discussed below and illustrated in Figure \ref{fig:Methodology}. A distinction is made between the strategic decisions of locating depots and assigning a number of vehicles to each depot (\textit{first stage}), and the reactive decisions of routing the vehicles for a specific realization of the demand (\textit{second stage}). By making this distinction between first and second stage, it is natural to decompose the problem and solve it as a two-stage stochastic program. The strategic or first-stage decisions are initially made without knowledge of demand realizations, then routing decisions are made for demand scenarios, and the locations are adjusted according to second stage results. For solving this problem we present a modification of a multi-cut \textit{L-Shaped} algorithm that establishes communication between the two stages. 

The first stage (green box in Figure \ref{fig:Methodology}) is modeled as a facility location and allocation problem in which the expected value of the second stage (i.e., future routing decisions) is approximated. This approximation is done via multiple evaluations of the second stage by querying a series of pre-generated set of scenarios (communication node 1 in Figure \ref{fig:Methodology}). In the second stage (blue box in Figure \ref{fig:Methodology}) a sub-problem corresponding to a Multi Depot Vehicle Routing Problem (MDVRP) must be solved (orange box in Figure \ref{fig:Methodology}) for each scenario. The purpose of each sub-problem is to evaluate the first stage solution under each future scenario and provide feedback of its quality.

To solve the MDVRP in each second-stage sub-problem, we propose a column generation scheme in which a set covering formulation selects a subset of routes that visit each demand node at minimum cost. The subset of routes is selected from routes generated in an iterative fashion by consulting a route generator (communication nodes 2 in Figure \ref{fig:Methodology}). We reformulate the problem of finding feasible routes for each depot and each scenario as an Elementary Shortest Path Problem with Resource Constraint and Replenishment (ESPPRC-R). These ESPPRC-R can be solved efficiently, after a graph transformation, with a slight modification to the pulse algorithm presented by \cite{Lozano2013} combining the ideas presented by \cite{Lozano2016} and \cite{Bolivar2014}.

Scenarios are generated in a three step process (red box in Figure \ref{fig:Methodology})which consists on i) understanding the statistical properties of the demand ii) simulating new observations exploiting those properties and iii) selecting worst-case observations as scenarios with a MIP. 

The main contribution of this research is an end-to-end methodology that addresses the complexity of coupling strategic decisions with the combinatorial nature of considering multiple routing problems in space and time, while acknowledging uncertainties on key parameters. Precisely, for each block of the solution methodology we propose key adaptations and upgrades that allow our methodology to tackle the VFSPRP-SC efficiently. First, to the best of our knowledge, we propose a new set of cuts for the L-Shaped algorithm. These cuts combine the ideas from the traditional combinatorial Benders cuts, that provide incentives to binary variables in the model to explore the solution space, with the lower bounds on the expected second stage values provided by the classic continuous Benders cuts. Second, we reformulate the problem of finding feasible routes for each depot and each scenario as a Shortest Path Problem with Resource Constraint and Replenishment (ESPPRC-R). The novelty of considering the ESPPRC with replenishment relies on the possibility of defining a set of routes that could not be defined with the simple ESPPRC. Third, we present several acceleration strategies for the proposed methodology: we introduce a set of valid inequalities with the purpose of saving up computational time at an early stage of the L-Shaped algorithm; we propose a strategy that invokes the route generators only when a promise of improvement exists; we implement an algorithm that allows for routes generated for a particular configuration of the demand to be mapped into another. Finally, we offer a tool that provides efficient data management, reporting, analytics, modeling and planning issues to improve logistics and distribution performance implemented in AIMMS and oriented to daily routing operation.

The paper is organized as follows. In \S \ref{s:StatAnalysis} we analyze the statistical properties of the data and place some context to it. In \S \ref{s:SolMethodology} we elaborate on the models and strategies that will comprise the solution methodology. In \S \ref{s:CaseStudy} we report and discuss results on how the model behaves under fluctuations of demand. Finally, in \S \ref{s:Conclusions} we provide conclusions and ideas for future work.

\newpage
\section{Statistical Analysis}\label{s:StatAnalysis}

Statistical analysis can give hindsight of the data on how to address the problem at hand and it is the first of the proactive steps towards business intelligence  \citep{davenport2007competing}. The purpose of this section is to explore statistical properties of the data that will support the development of the solution scheme for the selection of depots, determination of fleet size and vehicle routes.

\begin{figure}[h]
    \centering
    \includegraphics[width=0.4\textwidth]{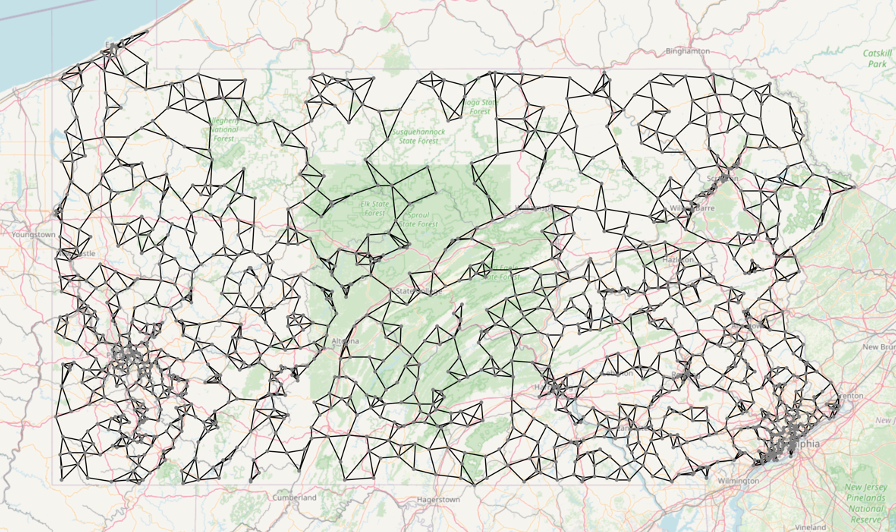}
    \caption{Network of locations and routes in the available data.}
    \label{fig:pennsylvania}
\end{figure}

Information about a set of locations is given in the form of 5-digit zip-codes which represent locations in the state of Pennsylvania, US. Moreover, these locations are connected with a set of routes between them. Figure \ref{fig:pennsylvania} presents the network that is created with the coordinates of the locations and their connections. At any given day, some locations are going to place orders, henceforth \textit{active clients}. The orders are going to vary on the number of units in them (order size).

\begin{figure}[h]
\centering
   \begin{subfigure}{0.4\textwidth}
     \centering
     \includegraphics[width=0.95\textwidth]{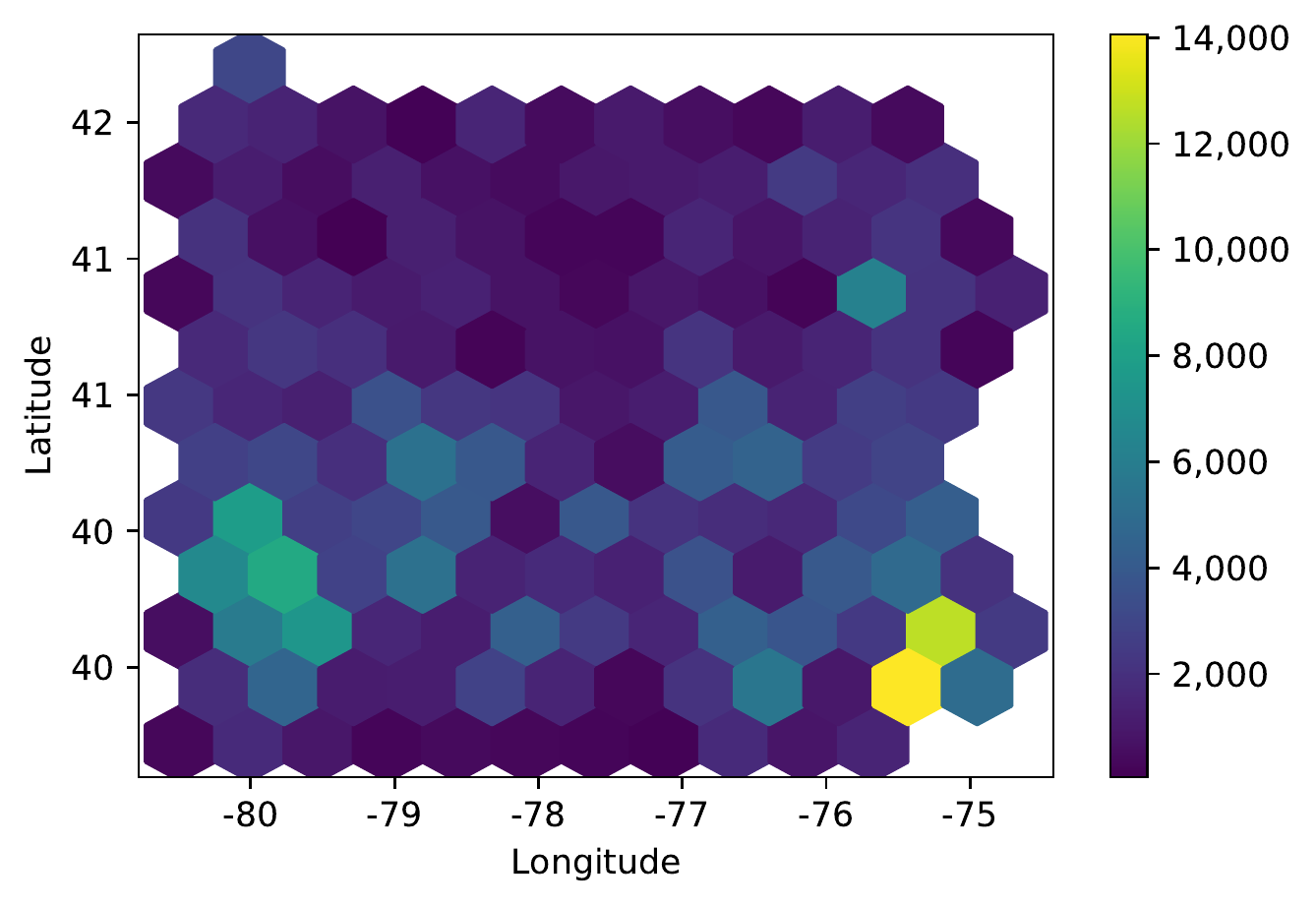}
     \caption{Spatial distribution of yearly demands for 2018.}
     \label{fig:Sdist2018}
   \end{subfigure}
   \begin{subfigure}{0.4\textwidth}
     \centering
     \includegraphics[width=0.95\textwidth]{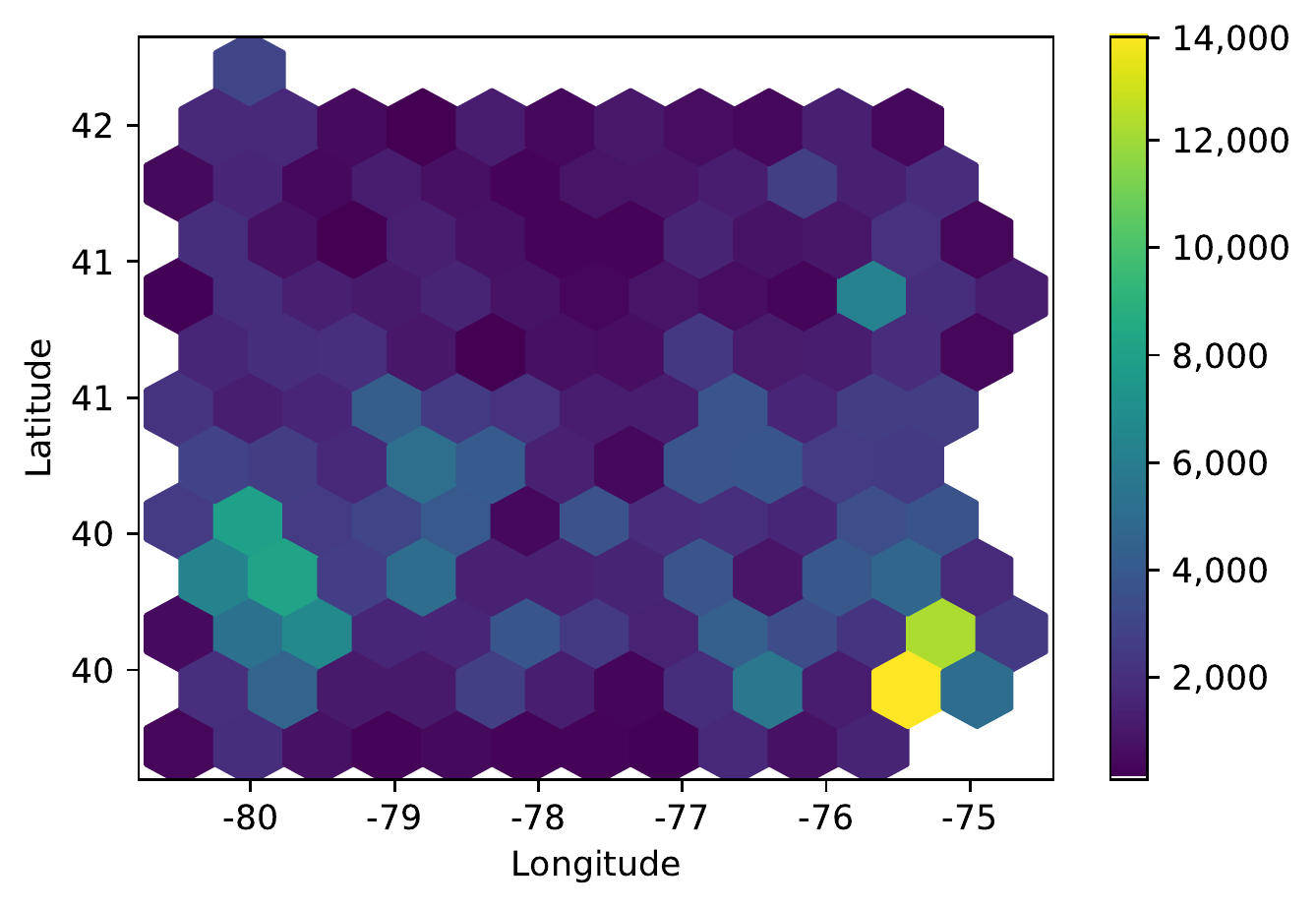}
     \caption{Spatial distribution of yearly demands for 2019.}
     \label{fig:Sdist2019} 
   \end{subfigure}
\caption{ Spatial and scalar distribution for the yearly demand of a customer.}\label{fig:Sdist}
\end{figure}

Figure \ref{fig:Sdist2018} and Figure \ref{fig:Sdist2019} show the aggregated yearly demand over the map with color signatures for 2018 and 2019, respectively. The color purple represents lower values and the color yellow represents the higher values. From Figure \ref{fig:Sdist} it can be noted that the behavior of the demand is similar in both years both spatially and in quantities. Some sectors of the network of locations that show larger yearly demand and these sectors correspond with the areas of the map with more locations. Particularly, the most active areas correspond the cities of Philadelphia and Pittsburgh.

\begin{figure}[h]
\centering
   \begin{subfigure}{0.4\textwidth}
     \centering
     \includegraphics[width=0.95\textwidth]{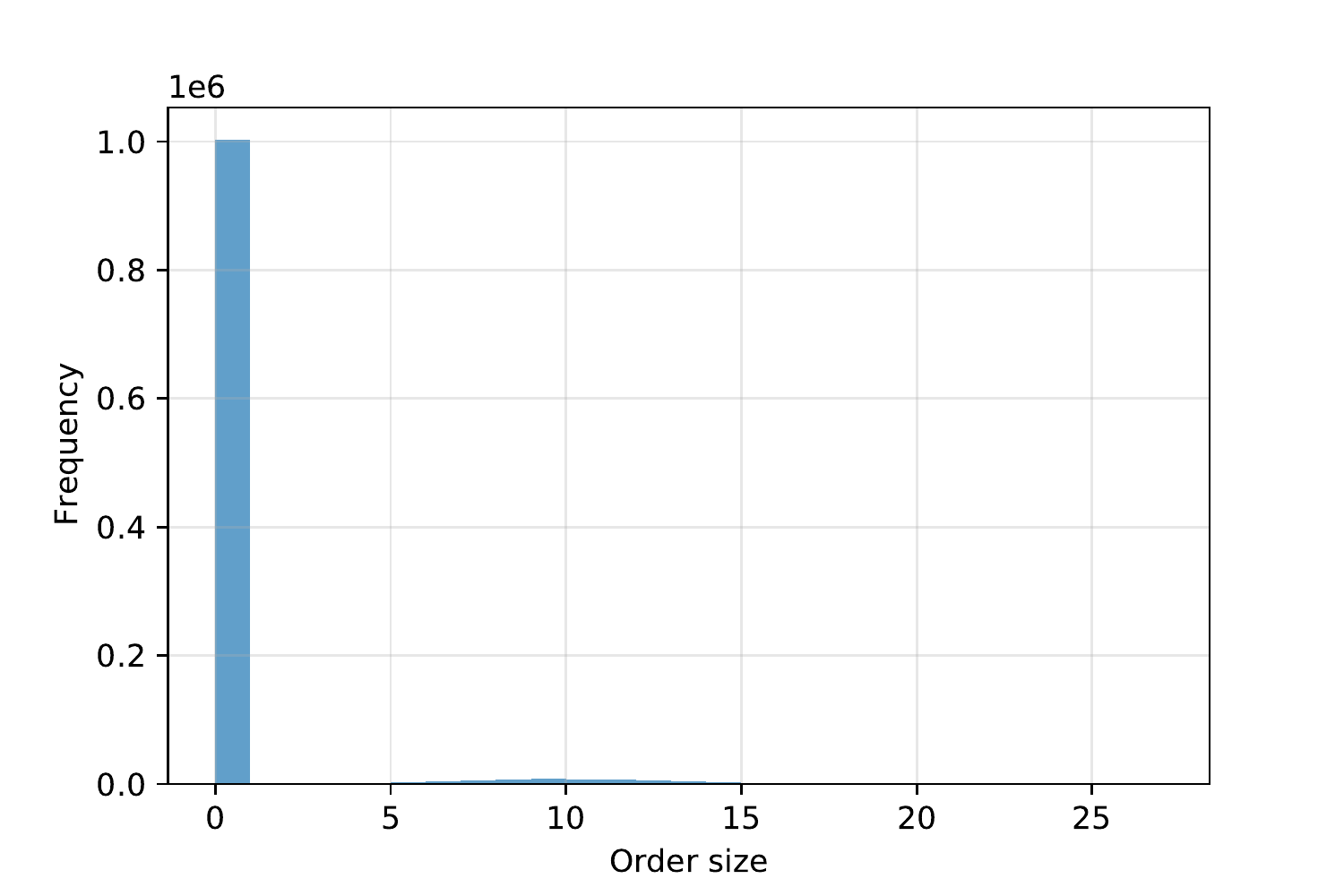}
     \caption{Distribution for the order size with zero counts.}
     \label{fig:PoisZero}
   \end{subfigure}
   \begin{subfigure}{0.4\textwidth}
     \centering
     \includegraphics[width=0.95\textwidth]{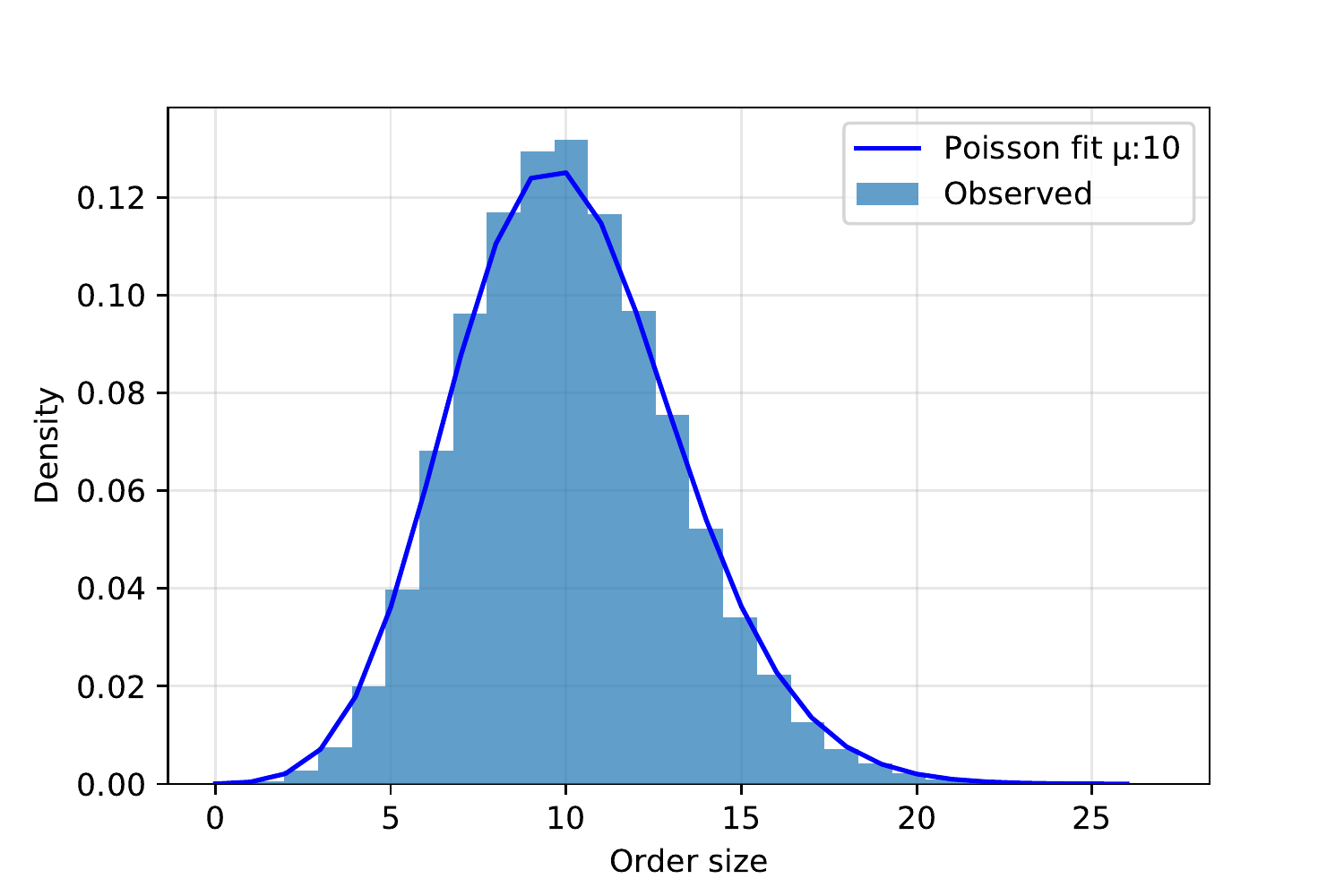}
    \caption{Distribution for the order size without zero counts.}
    \label{fig:PoisOrder}
   \end{subfigure}
\caption{Distribution for the order size.}\label{fig:PoisDist}
\end{figure}

\small
\normalsize
Figure \ref{fig:PoisZero} suggest the existence of patterns in the order size. Most of the time, the order size in any location is zero. A zero-inflated distribution can be better observed by splitting the distribution into two parts, i) a first zero generating process and ii) a process governed by other distribution that generates counts, some of which may be zero. Figure \ref{fig:PoisOrder} shows the distribution of the order size without the zero observations with a best-fit for a Poisson distribution with a mean of 10 units.

\small
\normalsize

\begin{figure}[h]
    \centering
    \includegraphics[width=0.5\textwidth]{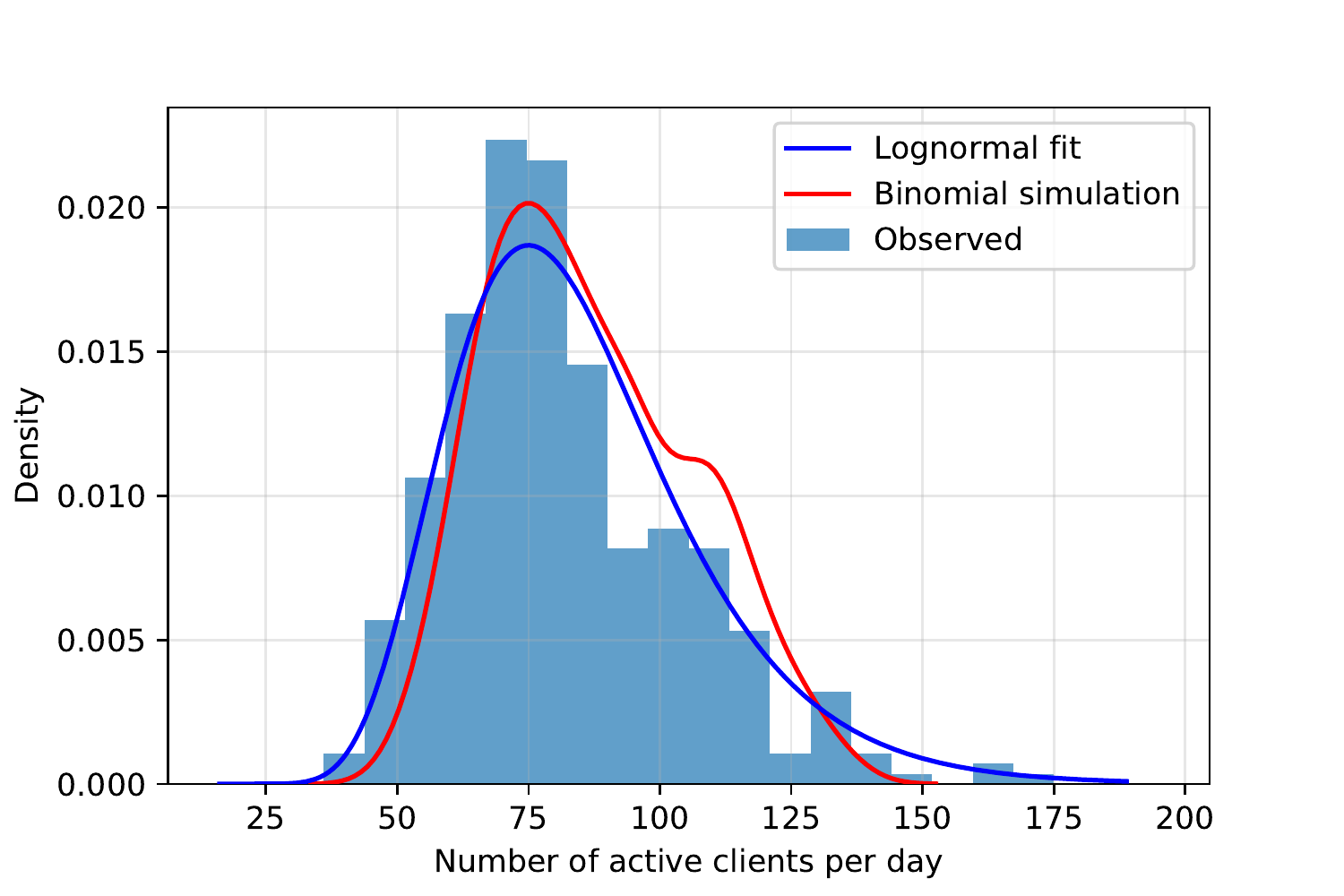}
    \caption{Distribution comparison for the daily number of active clients}
    \label{fig:active_clients}
\end{figure}

As for the zero generating process, two approaches were evaluated. The first consisting in simulating Bernoulli trials and the second consisting on best-fitting a non-discrete distribution. The efficacy of both approaches was compared using the distribution of active clients shown in Figure \ref{fig:active_clients} which is the complement of the distribution of the zero generating process. The distributions in both approaches were trained with 2018 data and tested over 2019 data. For the Bernoulli approach, a non-zero demand was considered a success and the probability of a success for a location in a month was calculated as the number of successes divided by the number of days in that month. In the second approach, a log normal distribution was the best-fit. The later seems more appropriate to adjust to the number of active clients. With this approach, the number of active clients can be simulated and by extension, the zero generating process.

\newpage
\section{Solution methodology}\label{s:SolMethodology}

In the proposed methodology, proactive decisions are made before the realization of uncertain events (demand), while reactive decisions are implemented afterwards. For the VFSPRP-SC, the first stage implies depot location and allocation of trucks, whereas the second stage addresses the routing of the previously allocated vehicles, given a realized demand.

Let $\G=(\N,\A)$ be a directed graph representing a transportation network, where $\N$ represents the set of nodes, and $\A \subseteq \N\times \N$ the set of arcs. Each arc $a\in \A$ has a distance $\delta_a$ (miles), a time $\tau_a$ (hours), and a cost $c_a$ (dollars) for a vehicle to move along the arc $a\in\A$. Let $\Hh \subseteq \N$ be the subset of nodes in which depots can be located, $h$ the number of depots that need to be opened and for $k\in \Hh$ let $f_k$ be the daily fixed cost of locating and maintaining a depot in that location. The maximum number of vehicles that can be allocated at any depot is denoted by $M$, and each vehicle allocated at any depot has a daily fixed cost of $b$. For a given realization of a random event $\omega\in\Omega$, the second-stage problem data, namely the demand locations, denoted by $\D(\omega)\subseteq\N$ and their order size, denoted by $\tilde{d}_i(\omega)$, $i\in\N$, become known. Piecing together the stochastic components of the second-stage data, we obtain the vector $\bm{\xi}^T(\omega)=\left(\D(\omega),(d^\omega_i)_{i\in\D(\omega)}\right)$. In this sense, a single random event $\omega$ (or state of the world) influences several random variables, namely, all components of $\bm{\xi}$. 

Let $x_k$ be a binary variable that takes the value of one if a depot is located at node $k\in \Hh$, and $y_k$ be the number of vehicles allocated at node $k\in\Hh$. For a given realization of the random event $\omega$ and a given first stage solution $(\bm{x},\bm{y})\in\{0,1\}^{|\Hh|}\times \mathbb{N}^{|\Hh|}$, the total cost of the MDVRP solved with the data $\bm{\xi}(\omega)$ will be denoted by $\mathcal{Q}(\bm{x},\bm{y},\bm{\xi}(\omega))$. The deterministic equivalent formulation for the VFSPRP-SC is then

\begin{mini!}|s|
    {}{\sum_{k\in\Hh} f_k x_k +\sum_{k\in\Hh}by_k +\mathbb{E}_\xi\left[\mathcal{Q}(\bm{x},\bm{y},\bm{\xi}(\omega))\right] \protect\label{eq:DEObj}}
    {\label{prob:DE}}{}
    \addConstraint{y_k}{\le Mx_k,}{\quad \forall k\in\Hh \protect\label{eq:DERel}}
    \addConstraint{x_k}{\in \{0,1\},}{\quad \forall k\in\Hh \protect\label{eq:DENatx}}
    \addConstraint{y_k}{\in \mathbb{N},}{\quad \forall k\in\Hh. \protect\label{eq:DENaty}}
\end{mini!}
\noindent Constraints (\ref{eq:DERel}) guarantee that vehicles can only be allocated at nodes where a depot is placed. In the objective function (\ref{eq:DEObj}) the first term accounts for the total cost of locating and maintaining the depots; the second term accounts for the total cost of buying and maintaining vehicles; and the third term captures the expected value of the \textit{second-stage} objective $\mathcal{Q}(\bm{x},\bm{y},\bm{\xi}(\omega))$ taken over all realizations of the random event $\omega\in\Omega$. The term $\mathcal{Q}(\bm{x},\bm{y},\bm{\xi}(\omega))$ is, itself, the solution to a MDVRP and needs to be solved for every $\omega\in\Omega$.

This section is organized as follows: In \S \ref{ss:Simulations} we present a simulation for scenarios of demand; in \S \ref{ss:LShaped} we present the L-Shaped algorithm used for solving the deterministic equivalent of the VFSPRP-SC; in \S \ref{ss:SecondStage} we present the proposed methodology for solving a MDVRP and the auxiliary problem that generates the routes; and finally in \S \ref{ss:AcelerationStrategies} we present several acceleration strategies for the algorithms presented in \S \ref{ss:LShaped} and \S \ref{ss:SecondStage}.

\subsection{Simulation}\label{ss:Simulations}

The purpose of this section is to present a method that makes computationally tractable the approximation of the nonlinear term in the objective function of problem (\ref{prob:DE}), namely 
\begin{equation}\label{eq:ExpSS}
\mathbb{E}_\xi\left[\mathcal{Q}(\bm{x},\bm{y},\bm{\xi}(\omega))\right].
\end{equation}
\noindent Since the realizations of vector $\bm{\xi}$ are numerous, solving the model for all of them would become intractable. For that matter we select a subset of scenarios to represent the variability of demand. For that we recall the descriptive models of \S \ref{s:StatAnalysis} to build a sample with simulation and define scenarios as shown in the bottom part of \ref{fig:Methodology}

One particularity of the problem is that solutions that satisfy the demand of all the active clients in any given day are preferred. If the scenarios were selected from the provided days of data then it would be not possible to foresee difficult scenarios (i.e., peaks in demand) that could affect strategic planning (i.e., the selection of depots). Difficult scenarios manifest as days with large demand, a high number of clients, scattered clients or a combination of the three conditions; and finding them is not trivial. Solving the resource allocation problem with worst-case scenarios improve the chances to achieve a good service level over time. With that in mind, we propose to solve the deterministic equivalent of the VFSPRP-SC using worse case scenarios selected with model \ref{Aps:scensel} from a simulated sample. 

\begin{figure}[h]
    \centering
    \includegraphics[width=0.5\textwidth]{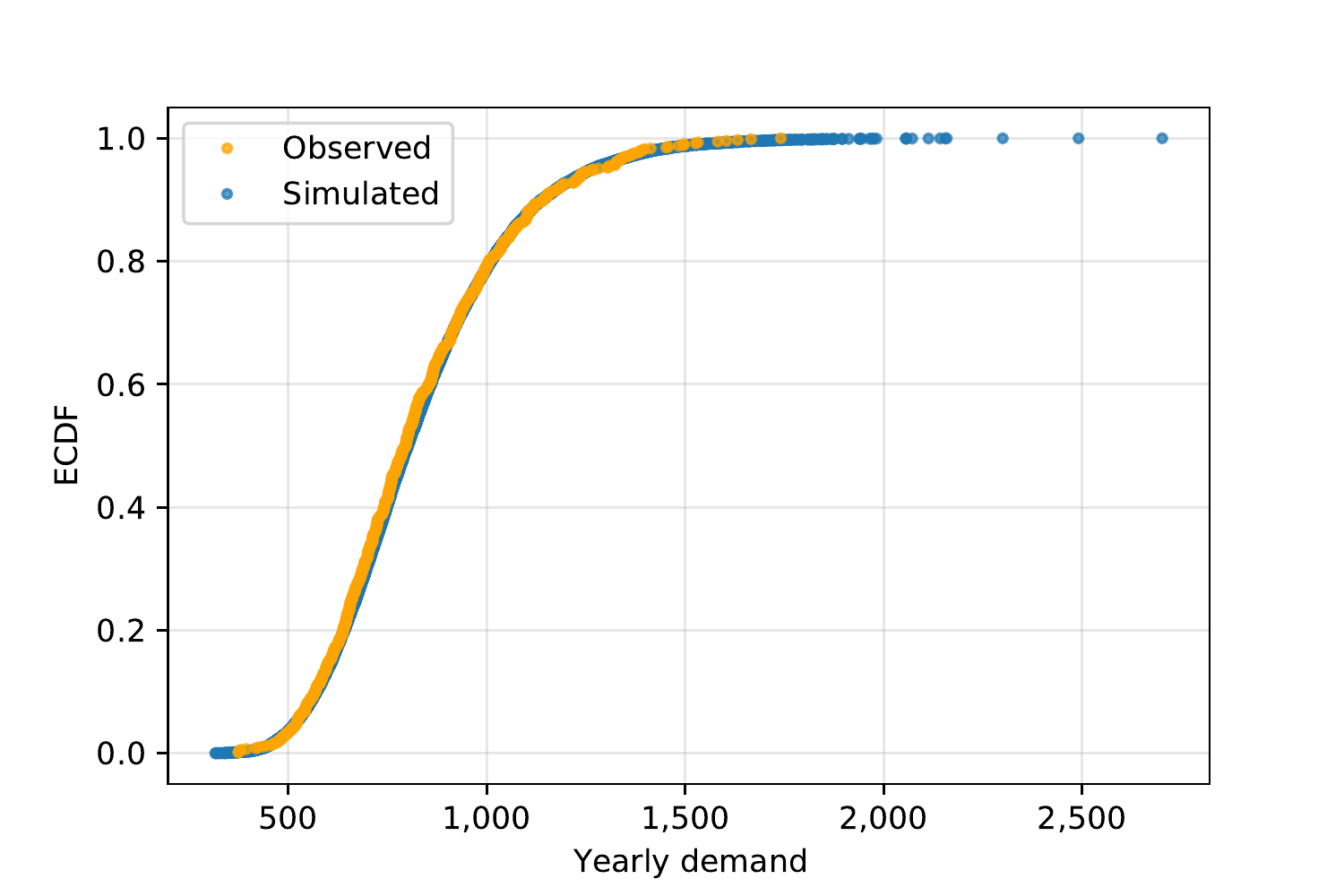}
    \caption{Empirical cumulative distribution function for the daily aggregated yearly demand of observed data and simulated data (10,000 samples).}
    \label{fig:ecdf}
\end{figure}

We simulate the demand by combining the adjusted Poisson and log normal distributions explored in \S \ref{s:StatAnalysis}. The log normal distribution gives the number of active clients in each observation. The active clients are going to be selected using their Bernoulli probabilities and the order size is assigned with the Poisson distribution. Figure \ref{fig:ecdf} shows that not only the simulated data (in blue) adjusts to the observed data (in orange) but also considers larger demands. 

Then, we determine ranges in the simulated cumulative distribution to become scenarios. An observation can be classified within the sample of a scenario if its demand fall in its range. Model \ref{Aps:scensel} selects from the pool of observations of each scenario the one with the highest number of active clients in a predefined group of geographical divisions to use as the most robust observation. The number of clients and the demand are highly correlated, therefore, one of the two is sufficient as the function to optimize. If the number of active clients is maximized for every geographical divisions, the result is an observation with scattered demand. We selected districts of Pennsylvania as the geographical divisions because they are based in population. We use the generated scenarios to train the first stage of the model as described below. 





\subsection{A Multicut \textit{L-Shaped} Algorithm for the VFSPRP-SC}\label{ss:LShaped}

With the results obtained in \S \ref{ss:Simulations}, we now formulate an estimation of \eqref{prob:DE} replacing the expected value of the second stage \eqref{eq:ExpSS} by
\begin{equation}
    \sum_{s\in \Scen}p_s\eta_s
\end{equation}
where each $\eta_s$ is a continuous variable that represents our estimation of $Q(\bm{x},\bm{y},\bm{\xi}(s))$\citep{birge2011}. This approximation will be called the \textit{master program} (MP) and is the following

\begin{mini!}|s|
    {}{\sum_{k\in\Hh} f_k x_k +\sum_{k\in\Hh}cy_k +\sum_{s\in \Scen}p_s\eta_s \protect\label{eq:BMPObj}}
    {\label{prob:BMP}}{\nu^*=}
    \addConstraint{y_k}{\le Mx_k,}{\quad \forall k\in\Hh \protect\label{eq:BMPRel}}
    \addConstraint{\eta_s}{\ge \sum_{i\in\D^s}\pi^{t,s}_i +\sum_{k\in\Hh}\lambda^{t,s}_ky_k - \eta^{t^*}_s\left(\sum_{k\in {O}^{t}}(1-x_k) +\sum_{k\not\in {O}^{t}}x_k \right)}{\quad \forall s\in\Scen,t\ge0 \protect\label{eq:BMPOpt}}
    \addConstraint{x_k}{\in \{0,1\},}{\quad \forall k\in\Hh \protect\label{eq:BMPNatx}}
    \addConstraint{y_k}{\in \mathbb{N},}{\quad \forall k\in\Hh \protect\label{eq:BMPNaty}}
    \addConstraint{\eta_s}{\ge 0,}{\quad \forall s\in\Scen. \protect\label{eq:BMPNatEta}}
\end{mini!}

Here constraints \eqref{eq:BMPOpt} are the actual approximation to the second stage value at each $s\in\Scen$, and will be added by consulting the second stage sub-problems (i.e., MDVRP's). The index $t=0,\cdots\ $ represents the number of evaluations of the second stage. Note that the term 
\begin{equation}\label{eq:DRRMPObj}
    \sum_{i\in\D^s}\pi^{s,t}_i +\sum_{k\in\Hh}\lambda^{s,t}_ky_k
\end{equation}
\noindent in constraint \eqref{eq:BMPOpt} corresponds to the classic continuous optimality cut for a problem in which the only variables that influence the second stage were the $y_k$ for $k\in\Hh$ \citep{birge2011}. Here $\pi^{s,t}_k$  and $\lambda^{s,t}_k$ represent the dual multipliers of the optimal solution of the for the $t$-th query  sub-problem associated with scenario $s\in\S$ (see \S \ref{ss:SecondStage}). Since the second stage is also influenced by the variables $x_k$ for $k\in\Hh$, we propose to subtract the term 
\begin{equation}\label{eq:BMPCombCut}
    \eta^{t^*}_s\left(\sum_{k\in {O}^{t}}(1-x_k) +\sum_{k\not\in {O}^{t}}x_k \right)
\end{equation}
\noindent to \eqref{eq:DRRMPObj} for obtaining the cut \eqref{eq:BMPOpt}, where $O^t\triangleq\{k\in\Hh: \hat{x}^t_i=1\}$, and $\eta^{t^*}_s=Q(\bm{\hat{x}}^t,\bm{\hat{y}}^t,\bm{\xi}(s))$.

Note that \eqref{eq:DRRMPObj} corresponds to the objective function of the dual program of the relaxed set covering formulation for the MDVRP (see \S \ref{ss:SecondStage}). Since the continuous relaxation is a linear program, the strong duality theorem holds, and for the optimal values of $\bm{\lambda}^{t,s}$ and $\bm{\pi}^{t,s}$, and $\bm{\hat{y}}$ fixed, \eqref{eq:DRRMPObj} equals to the optimal objective value of the relaxed set covering formulation for the MDVRP. Therefore, this cut acts as the combinatorial Benders cut, in the sense that unless some new combination for the binary $\bm{x}$ variables is tried, the second stage estimation for scenario $s\in\Scen$ needs to be grater than \eqref{eq:DRRMPObj}. However, if no new combination of the binary variables is tried, the proposed cut acts as a classical continuous Benders cut, using the information of the dual multipliers of the sub problems for estimating a lower bound on the second stage cost of scenario $s\in\Scen$ as function of the $\bm{y}$ variables.

In summary, cut \eqref{eq:BMPOpt} forces $\eta_s$, the second stage estimation for scenario $s\in\Scen$, to be grater than \eqref{eq:DRRMPObj}, the classic continuous Benders cut, unless some new combination of depot location is tried. 

In Appendix \ref{Ap:Algorithms}, Algorithm \ref{alg:LShaped} outlines the L-shaped algorithm followed to solve the VFSPRP-SC. This algorithm establishes a communication between the first and the second stage in which every possible solutions is consulted trough the $(\bm{\hat{x},\hat{y}})$ variables with the second stage subproblems, and feedback of the queried first stage solution is sent in form of the proposed \eqref{eq:BMPOpt} cuts, as outlined in Figure \ref{fig:Methodology}.

\subsection{On solving the Multi-depot Vehicle Routing Problem }\label{ss:SecondStage}

For solving each MDVRP we use the set covering approach proposed by \cite{Contardo2014}. For a given solution of the first stage $(\hat{\bm{x}}$, $\hat{\bm{y}})\in\{0,1\}^{|\Hh|}\times\mathbb{N}^{|\Hh|}$ and for a given realization of the random event $\omega \in\Omega$, let $\Rout(\omega)$ be the set of feasible routes for the second-stage problem data $\bm{\xi}(\omega)$. For depot $k\in \Hh$, we denote by $\Rout_k(\omega)$, the subset of routes staring and ending at $k$. 
For each $r\in\Rout(\omega)$, we let $z_r$ be a binary variable that equals to 1 if the route $r\in\Rout(\omega)$ is selected, and we denote by $c_r$ its cost, which equals to the sum of traveling costs along the arcs used by $r$, $\left(\text{i.e. }\sum_{a\in r}c(a)\right) $. In addition, we define the indicator function $e_r:\N\mapsto \{0,1\}$ that takes the value of one if the node $i\in\N$ is part of the route $r\in\Rout(\omega)$ and zero otherwise. The set covering formulation for the MDVRP is as follows

\begin{mini!}|s|
{}{\sum_{r\in\Rout(\omega)}c_r z_r \protect\label{eq:SPMfo}}
{\label{prob:SPM}}{\mathcal{Q}(\hat{\bm{x}},\hat{\bm{y}},\bm{\xi}(\omega))=}
\addConstraint{\sum_{r\in\Rout(\omega)}e_r(i) z_r}{\ge 1,\quad}{\forall i\in\mathcal{D}(\omega)\quad (\pi_i)\protect\label{eq:SPMset_part}}
\addConstraint{\sum_{r\in\Rout_k(\omega)}z_r}{\le \hat{y}_k,\quad}{\forall k\in\Hh \quad (\lambda_k) \protect\label{eq:SPMcap}}
\addConstraint{z_r}{\in\{0,1\},\quad}{\forall r\in\Rout(\omega). \protect \label{eq:SPMnatur}}
\end{mini!}

\noindent The objective function (\ref{eq:SPMfo}) aims to minimize the total traveling cost; constraints (\ref{eq:SPMset_part}) are the set covering constraints that ensure that each customer with a positive demand is visited at least once; and constraints (\ref{eq:SPMcap}) impose the first stage decisions, i.e., the depots opened, and the number of vehicles allocated at each depot. Constraints (\ref{eq:SPMnatur}) state the binary nature of the variables and by replacing them with 
\begin{equation}\label{eq:relaxed}
    0\le z_r\le 1,\quad\quad \forall r\in\Rout(\omega)
\end{equation}
we obtain a continuous relaxation of the problem. This relaxation allow us to obtain the dual variables of constraints \ref{eq:SPMset_part}, labeled as $\pi_i$ for $i\in \D(\omega)$, and of constraint \ref{eq:SPMcap}, labeled as $\lambda_k$ for $k\in\Hh$.

Note that the cardinality of the set of routes $\Rout(\omega)$ becomes prohibitive, and we cannot explicitly state all the variables of problem (\ref{prob:SPM}). In this case, the appealing idea of column generation (CG) to work only with a sufficient meaningful subset of routes becomes attractive. When formulating (\ref{prob:SPM}) and replacing $\Rout(\omega)$ with a subset of routes $\Rout'(\omega)\subseteq\Rout(\omega)$ we obtain the restricted master problem (RMP), and when considering the continuous relaxation of the RMP we obtain the restricted relaxed master problem (RRMP). The idea behind the CG algorithm is, like in the simplex method, to find in every iteration a promising variable to enter the basis, namely a route with negative reduced cost. The main difference between the CG and the Simplex method is that in the CG scheme these variables are not known (since the cardinality of the set $\Rout(\omega)$ is too big) and need to be generated. In this case, the reduced cost of a variable associated with the route $r\in\Rout_k(\omega)$ --the set of routes departing and ending at node $k\in\Hh$-- is given by

\begin{equation}\label{eq:reducedCost}
    \bar{c}_r=c_r -\lambda_k - \sum_{i\in\D(\omega)} a_r(i)\pi_i.
\end{equation}

In the Appendix \ref{Ap:Algorithms}, Algorithm \ref{alg:ColumnGeneration} outlines the specific CG followed to solve the MDVRP.  Given a specific first-stage solution $\hat{\bm{x}},\hat{\bm{y}}$ and second-stage problem data $\bm{\xi}(\omega)$ the procedure estimates the second-stage value $\mathcal{Q}(\hat{\bm{x}},\hat{\bm{y}},\bm{\xi}(\omega))$. In every iteration of the algorithm multiple auxiliary problems, one for each depot $k\in\hat{\Hh}$, use the dual multipliers of constraints \eqref{eq:SPMset_part} and \eqref{eq:SPMcap} to update the set of routes $\Rout'(\omega)$ from where the RMP selects, as outlined in Figure \ref{fig:Methodology}. This update is made with new and high quality routes, namely, routes with negative reduced cost. This procedure is repeated until the auxiliary problems can not find any more routes with negative reduced cost.

\subsubsection{Column generation Auxiliary problem (Route Generator)}\label{ss:RouteGenerator}

The purpose of this section is to present how given a depot $k\in\Hh$ and some dual vectors $\bm{\pi}$ and $\bm{\lambda}$ a route $r\in \Rout_k(\omega)$ with minimum reduced cost $\bar{c}_r$ is generated. First we must formally define the elements of $\Rout_k(\omega)$. 

Let's start by defining a new graph $\G^{sp}=(\N,\A^{sp})$, where the set of nodes is the same as the set of nodes of graph $\G$, namely $\N$, and $\A^{sp}=\left\{(v_i,v_j): v_i\in \N, v_j\in \N, v_i\ne v_j\right\}$, i.e., we consider a complete graph with the original set of nodes. For $a=(v_i,v_j)\in\A^{sp}$ let the $\delta^{sp}_a$ be the length of the shortest path from $v_i\in \N$ to $v_j\in \N$ in the original graph $\G$. Similarly we define the $c^{sp}_a$ and $\tau^{sp}_a$ as the cost and time of the shortest path from $v_i\in \N$ to $v_j\in \N$ over the original graph $\G$. A route is defined as a sequence of nodes $r=(v_0,v_1,\cdots,v_n)$, where $v_0=v_n$ and $v_i\in\N$ for $i=0,\cdots,n$, to which corresponds a sequence of arcs $(a_1,\cdots,a_n)$ where $a_i\in\A^{sp}$ for $i=1,\cdots,n$. At each node in the route the quantity $d_{v_i}(\omega)$ must be delivered. For a route departing from $k\in\Hh$ to be feasible, i.e., for a route $r$ to be in $\Rout_k(\omega)$, the following two conditions need to hold:

\begin{enumerate}
    \item Meet the maximum loading capacity per vehicle. This means that for each subsequence $(v_i,v_{i+1},\cdots,v_j)$ of the route $r$ such that $v_i=v_0=v_j$ and $v_k\ne v_0$ for $k=i+1,\cdots,j-1$, we must have that
    \begin{equation}
        \sum_{\ell=i}^j d_\ell(\omega) \le Q
    \end{equation}
    where $Q$ is the loading capacity per vehicle.
    \item Meet the time windows constraints for each node in the route. The vehicle reaches node $v_i$ at an hour 
    \begin{equation}
        t_i=t_0+\sum_{\ell=1}^i\tau^{sp}_{a_\ell}.
    \end{equation}
    Where $t_0$ is the time at which the vehicle departs from the depot. So the time windows constraint is met if for every node $v_i$ in the route $r$ if
    \begin{equation}
        w_1(v_i) \le h_i\le w_2(v_i)
    \end{equation}
    where $[w_1(v_i),w_2(v_i)]$ are the time windows in which node $v_i\in\N$ can  be served.
\end{enumerate}

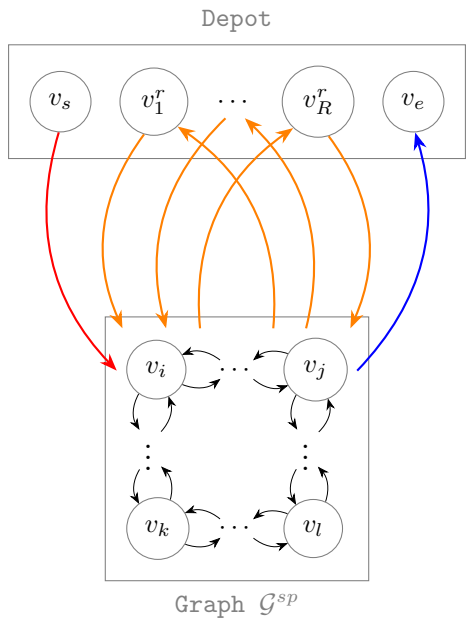
\begin{figure}
    \centering
    \begin{tikzpicture}[font=\ttfamily,
      mymatrix/.style={matrix of nodes, nodes=typetag, row sep=1em},
      mycontainer/.style={draw=gray, inner sep=1ex},
      typetag/.style={circle,draw=gray, inner sep=1ex, anchor=west},
      title/.style={draw=none, color=black, inner sep=0pt},myarrowR/.style={-Stealth,thick,bend right,black},
      myarrowL/.style={-Stealth,thick,bend left,black}
      ]
      \matrix[mymatrix] (mx1) {
        $v_s$ & |[title]|$\quad$ &
        $v^r_1$ &|[title]|$\quad$ & |[title]|$\cdots$ &|[title]|$\quad$ & $v^r_R$ &|[title]|$\quad$ &$v_e$\\
      };
       \matrix[mymatrix, below=4cm of mx1, matrix anchor=center] (mx2) {
        $v_i$ &|[title]|$\quad$ & |[title]|$\cdots$ &|[title]|$\quad$ & $v_j$\\
        |[title]|$\vdots$ & |[title]|$\quad$& &|[title]|$\quad$ & |[title]|$\vdots$\\
        $v_k$ &|[title]|$\quad$ & |[title]|$\cdots$ &|[title]|$\quad$ & $v_l$\\
      };
      \node[mycontainer, fit=(mx1)] {};
      \node[mycontainer, fit=(mx2)] {};
      \draw[myarrowR,red](mx1-1-1)edge(mx2);
      
      \draw[myarrowR,orange](mx1-1-3)edge(mx2);
      \draw[myarrowR,orange](mx2)edge(mx1-1-3);
      \draw[myarrowR,orange](mx1-1-5)edge(mx2);
      \draw[myarrowR,orange](mx2)edge(mx1-1-5);
      \draw[myarrowL,orange](mx1-1-7)edge(mx2);
      \draw[myarrowL,orange](mx2)edge(mx1-1-7);
      \draw[myarrowR,blue](mx2)edge(mx1-1-9);
      
      \draw[myarrowR,very thin](mx2-1-1)edge(mx2-1-3);
      \draw[myarrowR,very thin](mx2-1-3)edge(mx2-1-1);
      \draw[myarrowR,very thin](mx2-1-3)edge(mx2-1-5);
      \draw[myarrowR,very thin](mx2-1-5)edge(mx2-1-3);
      
      \draw[myarrowR,very thin](mx2-3-1)edge(mx2-3-3);
      \draw[myarrowR,very thin](mx2-3-3)edge(mx2-3-1);
      \draw[myarrowR,very thin](mx2-3-3)edge(mx2-3-5);
      \draw[myarrowR,very thin](mx2-3-5)edge(mx2-3-3);
      
      \draw[myarrowR,very thin](mx2-1-1)edge(mx2-2-1);
      \draw[myarrowR,very thin](mx2-2-1)edge(mx2-1-1);
      \draw[myarrowR,very thin](mx2-2-1)edge(mx2-3-1);
      \draw[myarrowR,very thin](mx2-3-1)edge(mx2-2-1);

      \draw[myarrowR,very thin](mx2-1-5)edge(mx2-2-5);
      \draw[myarrowR,very thin](mx2-2-5)edge(mx2-1-5);
      \draw[myarrowR,very thin](mx2-2-5)edge(mx2-3-5);
      \draw[myarrowR,very thin](mx2-3-5)edge(mx2-2-5);
      \node[gray,above=2mm of mx1]{$\text{Depot}$};
      \node[gray,below=2mm of mx2]{Graph $\G^{sp}$};
    \end{tikzpicture}
    \caption{Graphic representation of $\G^a$. Red arcs connect the start node $v_s$ with $\G^{sp}$; orange arcs connect the replenishment nodes $\{v^r_i\}_{i=1}^R$ with $\G^{sp}$; and blue arcs connect $\G^{sp}$ with the end node.}
    \label{fig:auxG}
\end{figure}

The problem of finding a route in $\Rout_k(\omega)$ is reformulated as an Elementary Shortest Path Problem with Resource Constraints and Replenishment (ESPPRC-R). This problem consists on finding an elementary path (i.e., ordered sequence of non repeating nodes) from a source node $v_s$ to an end node $v_e$ that minimizes the reduced cost subject to vehicle capacity and time windows constraints. In addition the ESPPRC-R considers replenishment nodes that reset the capacity of the vehicle if visited \citep{Lozano2016,Bolivar2014}. 
For the solutions of the ESPPRC-R to be feasible routes, we must solve the problem over an auxiliary Graph $\G^a=(\N^a,\A^a)$ (see Figure \ref{fig:auxG}). To the set of nodes of $\G^a$ we add $R+2$ auxiliary auxiliary nodes that represent the depot. The auxiliary nodes $v_s$ and $v_e$ represent the first and last node of the route, respectively; and the auxiliary nodes $v_1^r,\cdots,v^r_R$ are the replenishment or reloading nodes, and represent the visits to the depot in the middle of a route. The number of replenishment nodes $R$ is defined as
$$R=\ceil*{\frac{\sum_{i\in\D(\omega)}d_i^\omega}{Q}},$$
the maximum number of replenishment nodes possibly needed. Consequently, $\N^a$ is defined as $\N\cup\{v_s,v_t\}\cup\{v^r_1,\cdots,v^r_R\}$. To connect these nodes we add to the set of arcs $\A^{sp}$, arcs from $v_s$ to each node in $\N$, namely $\A^s=\{v_s\}\times\N$ (painted in red in figure \ref{fig:auxG}); arcs from every node $v_i\in\N$ to the end node $v_t$, namely $\A^t=\N\times\{v_t\}$ (painted in blue in figure \ref{fig:auxG}); and arcs between every replenishment node $v^r_k$ for $k=1,\cdots,R$ and every node $v_i\in\N$, namely $\A^r=\left(\{v^r_i\}_{i=1}^R\times \N\right) \cup\left(\N\times \{v^r_i\}_{i=1}^R\right)$ (painted in orange in figure \ref{fig:auxG}). Consequently, $\A^a$
 is defined as $\A^a=\A^{sp}\cup \A^s\cup\A^t\cup\A^r.$

For solving the ESPPRC-R we use a modification of the pulse algorithm presented by \cite{Lozano2013}. By combining the ideas presented by \cite{Lozano2016} and \cite{Bolivar2014}, we obtain an efficient way to solve the ESPPRC-R, and therefore finding feasible routes for the MDVRP. The difficulty of solving this ESPPRC-R arises when solving such problems under a CG scheme. Since our objective is to find a route with minimum reduced cost, the weight of each arc is $\bar{c}_{ij}=c_{ij}-\pi_j$ for $(i,j)\in\A^{sp}$. When subtracting $\pi_j$ we sometimes have negative weights that make difficult to solve any shortest path problem. \cite{Lozano2016} present an efficient way to solve the ESPPRC under CG schemes by adapting the versatile pulse algorithm. The difference between their approach and the one we present here is that we consider valid for the solution of the ESPPRC to reload its capacity multiple times at the depot. For allowing this, we use the ideas presented in \citep{Bolivar2014} to solve ESPPRC-R, by solving a constraint shortest path problem with replenishment from node $v_s\in\N^a$ to node $v_t\in\N^a$ over the graph $\G^a$.


\subsection{Acceleration strategies}\label{ss:AcelerationStrategies}

In this section we propose several acceleration strategies to save up time in an early stages of the Algorithm \ref{alg:LShaped}. Note that in a in an initial stage, Algorithm \ref{alg:LShaped} hast not much information about the problem that wants to be solved, i.e., its approximation of the expected value of the second stage \eqref{eq:ExpSS} is no good, resulting in \textit{greedy} and \textit{shortsighted} first stage solutions that do not take into account the second stage. For example, in the first iteration of algorithm \ref{alg:LShaped} a valid first stage solution could be to allocate zero trucks in every depot, and estimate the second stage cost as zero for all $s\in\Scen$; clearly, this first stage solution would have very high second stage cost, and in further iterations, Algorithm \ref{alg:LShaped} would be forced to change that solution. The following strategies aim to tackle this problem and to take advantage of this poor first stage solutions to save up computational time in an early stage of Algorithm \ref{alg:LShaped}.

\subsubsection{Valid inequality for the master problem}\label{sss:AcelMp}

Since in general the VFSPRP-SC is looking to cover some demand in the network, and aims to be able to reach all demand nodes for any possible demand scenario, we force the MP to be able to do so. Let $V_i$ be the set of potential depots that could reach the node $i\in \N$ in less than a maximum time radius $T$, namely $V_i:=\{k\in\Hh:\tau^{sp}_{(k,i)}\le T \}$. So we add the binary variable $w_i$ for $i\in\N$ that takes the value of one if the node $i\in\N$ is reached by some depot in the network, and zero otherwise. For the proper activation of $w_i$ for $i\in\N$ we add constraint \eqref{eq:actwi}
\begin{equation}\label{eq:actwi}
        w_i\le \sum_{k\in V_i}x_{k},\quad \forall i\in\N,
\end{equation}
\noindent and add the term 
\begin{equation}\label{eq:wiPenalty}
    \sum_{i\in\N}\rho_i (1-w_i)
\end{equation}Sí
\noindent to the objective function, where $\rho_i$ is a penalty for not reaching node $i\in\N$. Depending on the values of $\rho_i$, the MP is forced to reach all nodes in the network or not.

\subsubsection{Route generators activation paradigm}

This route generators activation paradigm guides our main optimization procedure. We noticed that it is not necessary to generate new routes in every iteration of the main procedure, so we only generate routes when it is worth it. If for each scenario $s\in \Scen$, the set of routes $\Rout(s)$ has already high quality routes, a good first approximation on $\mathcal{Q}(\hat{\bm{x}},\hat{\bm{y}},\bm{\xi}(s))$ can be obtained by solving Model \eqref{prob:SPM} without calling the route generators. In this sense, we do need high quality routes in every moment of Algorithm \ref{alg:LShaped} but these do not need to be generated for low quality first stage solutions. To purposefully generate these routes, we propose the following ideas: first in every iteration of Algorithm \ref{alg:LShaped} we keep track of the set of routes generated in previous iterations for each $s\in\Scen$, so the routes generated for different first stage solutions are saved, and recycled for other similar first stage solutions; second, since at the beginning of the algorithm the set of routes is empty, we use the route generators as a preprocess to fill this set of routes with initial high quality routes for each scenario; and third, in order to keep generating high quality routes, we generate routes only for those first stage solutions that promise improvement, i.e. $(\hat{\bm{x}},\hat{\bm{y}})$ that update the primal bound of Algorithm \ref{alg:LShaped}.

\subsubsection{Petal Recycler}\label{sss:AcelSP}
Since in every iteration of algorithm \ref{alg:LShaped} we need to solve multiple times a MDVRP with different problem data, namely $\bm{\xi}(s)$ for $s\in\Scen$, we propose a way of recycling the routes from one second-stage problem data into another. This allows us to speed up the CG (algorithm \ref{alg:ColumnGeneration}) by staring with an initial high quality set of routes $\Rout^{s'}$. Our method maps the set of routes generated for a specific second-stage problem data $\bm{\xi}(s)$, namely routes in $r\in\Rout^s$ to feasible routes for other second-stage problem data $\bm{\xi}(s')$, namely $\Rout^{s'}$. 

\begin{figure}[htb]
    \centering
    \begin{tikzpicture}[depot/.style={draw=black,thin,fill=blue,fill opacity = 0.1,text opacity=1,inner sep=0pt,minimum size=10pt,font=\scriptsize},
    nr1/.style={circle,draw=black,thin,fill=green,fill opacity = 0.3,text opacity=1,inner sep=0pt,minimum size=5pt,font=\scriptsize},
    nr2/.style={circle,draw=black,thin,fill=orange,fill opacity = 0.3,text opacity=1,inner sep=0pt,minimum size=5pt,font=\scriptsize},
    nleg/.style={circle,draw=white,thin,fill=white,fill opacity = 0.3,text opacity=1,inner sep=0pt,minimum size=5pt,font=\scriptsize},
    arcr1/.style={-Stealth,thin},
    arcr2/.style={-Stealth,thin,orange},
    arcr3/.style={-Stealth,dashed,very thin},
    myarrow1/.style={-Stealth,thick,bend right,red},
    myarrow2/.style={-Stealth,thick,bend right,orange},
    myarrow3/.style={-Stealth,thick,white}]

    \node[depot](0) at (0,0){$h$};
	
    \node[nr1](1) at (-2,1){};
	\node[nr1](2) at (-3,2){};
	\node[nr1](3) at (-1,4){};
	\node[nr1](4) at (0.3,2){};
	\node[nr1](5) at (2,1){};
    
    \draw[arcr1](0)edge(1);
    \draw[arcr1](1)edge(2);
    \draw[arcr1](2)edge(3);
    \draw[arcr1](3)edge(4);
    \draw[arcr1](4)edge(5);
    \draw[arcr1](5)edge(0);
    
    \node[nr2](21) at (-2.5,0.5){};
	\node[nr2](22) at (-2.1,2.1){};
	\node[nr2](23) at (-0.8,3){};
	\node[nr2](24) at (1.1,2.2){};
	\node[nr2](25) at (1,1){};
    
    \draw[arcr2](0)edge(21);
    \draw[arcr2](21)edge(22);
    \draw[arcr2](22)edge(23);
    \draw[arcr2](23)edge(24);
    \draw[arcr2](24)edge(25);
    \draw[arcr2](25)edge(0);
    
    \draw[arcr3](1)edge(21);
    \draw[arcr3](2)edge(22);
    \draw[arcr3](3)edge(23);
    \draw[arcr3](4)edge(24);
    \draw[arcr3](5)edge(25);

    \node[nr1](l1) at (1.3,4){};
    \node[nr1](l2) at (2.3,4){};
    \node[nleg](l3) at (1.8,4.4){$r\in\Rout^s_{h}$};
    \draw[arcr1](l1)edge(l2);
    
    \node[nr2](l21) at (1.3,3.1){};
    \node[nr2](l22) at (2.3,3.1){};
    \node[nleg](l23) at (1.8,3.5){$r'\in\Rout^{s'}_{h}$};
    \draw[arcr2](l21)edge(l22);
    
    \end{tikzpicture}

    \caption{Petal recycler intuition. The green route is mapped to the orange route.}
    \label{fig:PetalRecycler}
\end{figure}
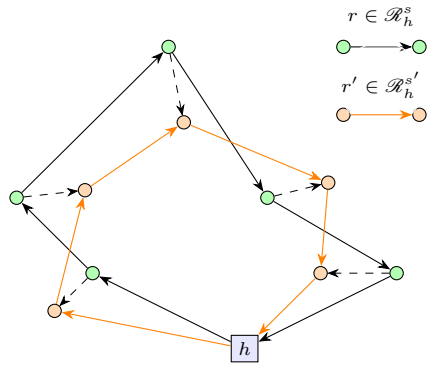

The intuition behind the methodology is very simple, given a route $r\in\Rout^s_k$, generated for the second-stage data $\bm{\xi}(s)$ and for a depot $k\in\N$, it seeks to replace every node in the route for a demand node of another second stage problem data $\bm{\xi}(s')$ (see Figure \ref{fig:PetalRecycler}). Clearly every replacement has to be made very carefully. First we want the new route to be as good as possible in terms of cost, and second, we want the new route to be feasible, i.e., to satisfy conditions 1 and 2 of \S \ref{ss:RouteGenerator}. In order to measure the cost of replacing the $i$-th node of route $r\in\Rout^s_k$, denoted by $v_{(i)}$, for a demand node $v_i\in\D(s')$ we define
$$\bar{c}(v_{(i-1)},v_{(i)},v_k)=\left(c^{sp}(v_{(i-1)},v_k)+c^{sp}(v_k,v_{(i+1)})\right)-\left(c^{sp}(v_{(i-1)},v_{(i)})+c^{sp}(v_{(i)},v_{(i+1)}) \right).$$

Note that $\bar{c}(v_{(i-1)},v_{(i)},v_k)$ gives us exactly the over cost or saving made by replacing node $v_{(i)}\in r$  for node $v_k\in\D(s)$, so if we replace nodes with low $\bar{c}(v_{(i-1)},v_{(i)},v_k)$ we guarantee that the quality of the route is kept. In order to keep the route feasible we only need to check before replacing the node, that the load and time windows constraints hold, as shown in algorithm \ref{alg:PetalRecycler}.

\section{Case study} \label{s:CaseStudy}

In this section we present the results for the competition instances. With the purpose of offering a more complete solution and alternatives to the decision maker, we use the flexibility of our methodology and provide a solution that considers the whole distribution, henceforth the expected value solution, and a solution that considers only the right tail of the demand distribution, henceforth the conservative solution. The expected value solution is obtained by solving model \eqref{prob:BMP} with a set of scenarios that addresses the average behavior of the demand, and does not consider \textit{right-tailed} demand scenarios. In contrast, the conservative solution is obtained by solving model \eqref{prob:BMP} with a set of scenarios that seeks a robust solution, capable of responding to extreme demand scenarios, even worse than the observed in the historical data. Both sets of scenarios were generated with the methodology outlined in \S \ref{ss:Simulations}. In order to compare these two solutions, a distinction is made between the deterministic part of the total daily operation costs, regarding to the strategic decisions; and the stochastic part, regarding to the routing costs. The estimation of this stochastic part, the routing costs, is done by solving for each day of the years 2018 and 2019 a MDVRP, and keeping track of its routing performance.

The rest of this section is organized as follows. \S \ref{ss:results} shows the strategic decisions suggested by our methodology; \S \ref{ss:SolEvaluation} explains how a proper description of the daily routing costs is made, and uses this description to fairly compare the two proposed solutions; and \S \ref{ss:compTimes} reports our computational times and shows the impact of our acceleration strategies. The, the computational experiments presented in this section were coded in AIMMS with a module in Python and a module in Java. All experiments were executed on an Intel core i7-7500U CPU@2.70GHz with 8GB of RAM.

\subsection{Strategic decisions}\label{ss:results}
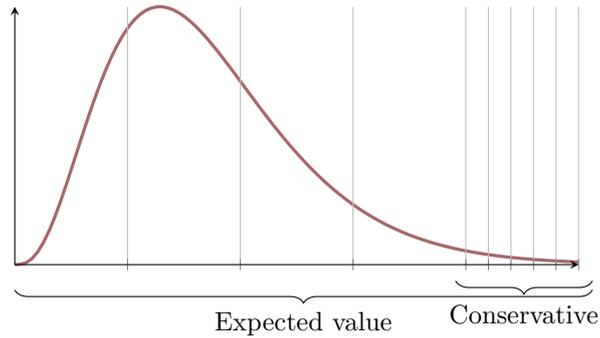
\begin{figure}[h]
            \centering
            \begin{tikzpicture}[
            declare function={gamma(\z)=
            2.506628274631*sqrt(1/\z)+ 0.20888568*(1/\z)^(1.5)+ 0.00870357*(1/\z)^(2.5)- (174.2106599*(1/\z)^(3.5))/25920- (715.6423511*(1/\z)^(4.5))/1244160)*exp((-ln(1/\z)-1)*\z;},
            declare function={gammapdf(\x,\k,\theta) = 1/(\theta^\k)*1/(gamma(\k))*\x^(\k-1)*exp(-\x/\theta);},
            declare function={f(\x) = 0;}
            ]
        
                \begin{axis}[
                    no markers, domain=0:9, samples=100,
                    axis lines=left, xlabel=$ $, ylabel=$ $,
                    every axis y label/.style={at=(current axis.above origin),anchor=east},
                    every axis x label/.style={at=(current axis.right of origin),anchor=north},
                    height=5cm, width=9cm,
                    xtick={7,14,21,28,29.4,30.8,32.2,33.6,35}, ytick=\empty,
                    xticklabels={},
                    enlargelimits=false, clip=false, axis on top,
                    grid = major
                    ]
                    \addplot [very thick,cyan!20!purple,domain=0:35] {gammapdf(x,4,3)};
                    
                \end{axis}
            \draw [decorate,decoration={brace,amplitude=5pt,mirror,raise=4ex}] (5.8,0.39) -- (7.6,0.39) node[midway,yshift=-3em]{Conservative};
            \draw [decorate,decoration={brace,amplitude=5pt,mirror,raise=4ex}] (0,0.27) -- (7.6,0.27) node[midway,yshift=-3em]{Expected value};
            \end{tikzpicture}
            \caption{Expected value solution vs. Conservative solution}
            \label{fig:RvAscenarios}
        \end{figure}

Figure \ref{fig:RvAscenarios} shows the percentiles for which the expected value and the conservative set of scenarios were selected. For the expected value solution, percentile in ranges from 20\% to 100\% were selected for simulating and choosing the scenarios. The purpose of the conservative set of scenarios is to provide robust strategic solution by considering rare events. For this set of scenarios the percentiles used for the selection were in the range of 90\% and 100\%.

\begin{figure}
    \centering
    \subfloat[Expected value solution \label{fig:solMapEV}]{\includegraphics[width=0.49\textwidth]{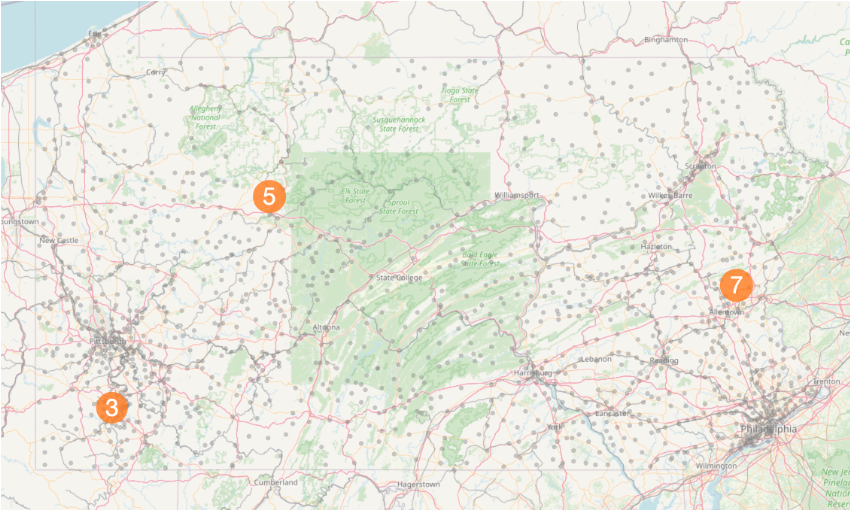}}
\subfloat[Conservative solution \label{fig:solMapC}]{\includegraphics[width=0.49\textwidth]{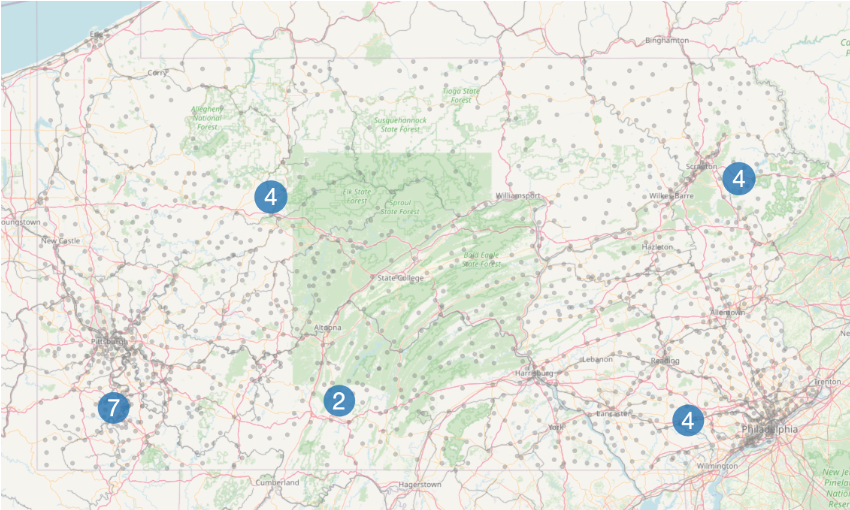}}
\\
    \caption{Spatial distribution of the results}
    \label{fig:solMap}
\end{figure}

Figure \ref{fig:solMap} outlines the spatial location of the selected depots for both set of scenarios. The depot-fleet size configuration obtained with the expected value scenarios opens 3 depots and allocates 15 vehicles in total, whereas the conservative solution opens 5 depots and allocates 21 vehicles in total. In both solutions we can see depots located on the outskirts of the main demand concentrations, this, probably aiming to reduce the transportation costs without incurring high depot maintaining costs within the urban areas. 

The expected value solution has daily depot maintenance costs of \$523 dollars, and vehicle maintenance costs of \$450 dollars, for a total daily fixed costs of \$973 dollars, whereas the conservative solution incur \$909 and \$630 dollars respectively for a total of \$1,539 dollars total daily fixed costs. As expected, the daily depot and vehicle maintenance costs of the expected value solution are almost 37\% cheaper than in the conservative case, this implies that the daily cost of being prepared for the worst-case scenarios is of \$566 dollars. This does not necessarily mean that the expected value solution is better, before jumping into conclusions we must first consider the daily routing costs of each of the solutions.

\subsection{Ex-post evaluation}\label{ss:SolEvaluation}

In order to accurately estimate the total daily costs, we must first understand how the routing costs of each of the solutions behave. For this task we performed the following ex-post evaluation. Since the routing costs depend on the stochastic demand, these are also a random variables. With this in mind, the only fair way of comparing the total daily expected cost of the operation for each of the solutions is to evaluate both solutions under the same set of demand scenarios, one different than the set of scenarios for which they were solved for. We could think the set of scenarios used in the optimization process as our \textit{training sets} and the set of scenarios in which we evaluate the solutions as our \textit{tests sets}. So, for evaluating the solutions, we use the historical data and we solved a MDVRP for each of the observed demand of the years 2018 and 2019. In \S \ref{sss:routingExpecctedCostst} we present the experiments performed for describing the random variables of the routing costs for each of the solution and we compare the daily routing performance. In \S \ref{sss:dailyExpecctedCostst} we compare the total daily costs of both solutions, considering the depot and fleet maintenance costs, and the routing costs. 

\subsubsection{Routing costs}\label{sss:routingExpecctedCostst}
Figure \ref{fig:routingCosts} and Table \ref{tab:routingCosts} report the summary of the ex-post evaluation. We use each MDVRP solved as a data point for describing the empirical probability density function of the routing costs of each of the solutions (Figure \ref{fig:routingCosts}). In Table \ref{tab:routingCosts} we present several descriptive statistics for both solutions. Column 1 reports the the descriptive statistics; column 2 and 3 report the value of each descriptive statistic for the Expected value solution and the conservative solution respectively.

\begin{table}[h]
    \begin{minipage}{0.5\linewidth}
		\centering
		\includegraphics[width=1\textwidth]{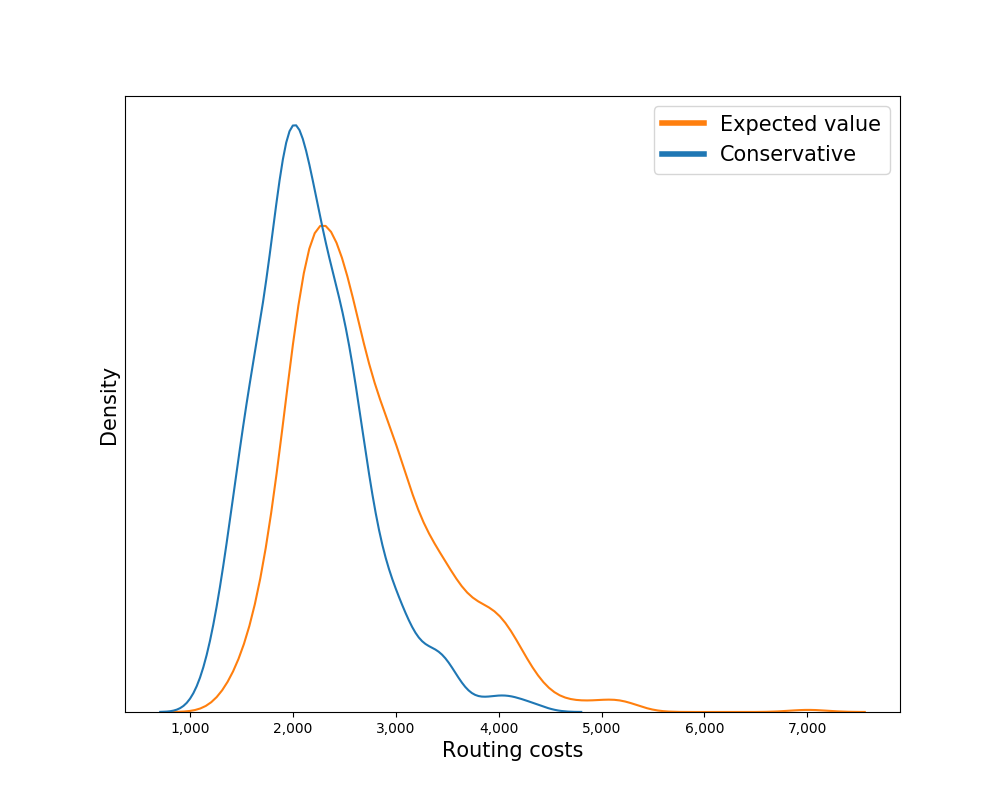}
		\captionof{figure}{Empirical probability density functions of the routing costs for both solutions. }
		\label{fig:routingCosts}
	\end{minipage}
    \begin{minipage}{0.5\linewidth}
        \caption{Descriptive statistics for the daily routing costs for both, the expected value and the conservative solutions}
		\label{tab:routingCosts}
		\centering\small
		\begin{tabular}{lrr}
            \hline
                               & \multicolumn{1}{c}{Exp. value} & \multicolumn{1}{c}{Conservative} \\ \hline
            Mean               & 2,693.21                                     & 2,194.64                                   \\
            Std. Error     & 26.31                                       & 19.91                                     \\
            Median             & 2,538.14                                     & 2,118.45                                   \\
            Std. Dev. & 710.49                                      & 537.9                                     \\
            Kurtosis           & 2.21                                        & 1.26                                      \\
            Skewness           & 1.14                                        & 0.88                                      \\
            Minimum            & 1,385.99                                     & 1,124.01                                   \\
            Maximum            & 7,009.36                                     & 4,390.77                                   \\ \hline
        \end{tabular}
    \end{minipage}
\end{table}

As Table \ref{tab:routingCosts} shows, the expected routing costs of the conservative solution is around 500 dollars cheaper than the expected routing costs of the expected value solution. When considering the standard deviation, we confirm that the risk taken by the conservative solution is lower than the one taken by the expected value solution, this can also be seen when comparing the worst case for both solutions, the maximum cost that the expected solution incur is around \$26 hundred  dollars garter than the conservative case.

\begin{figure}[h]
     \centering
     \includegraphics[width=0.7\textwidth]{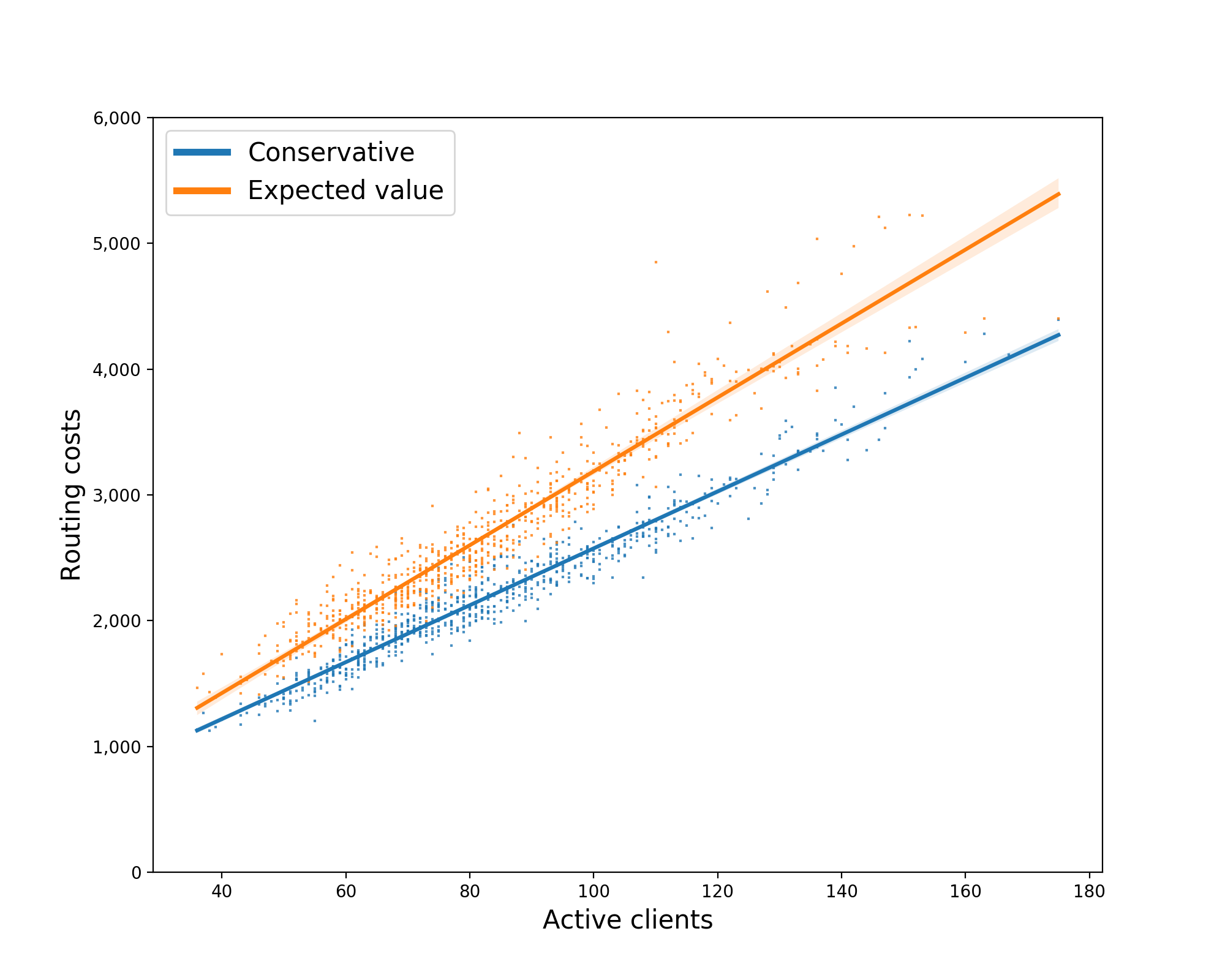}
     \caption{Regression plot for the routing costs explained by the active clients.}
     \label{fig:RegPlotRCvsAC}
 \end{figure}

 We noticed that the daily routing costs could be explained by the number of active clients (see Figure \ref{fig:RegPlotRCvsAC}), and we can estimate the average marginal cost of serving an extra client for each of the proposed solution. With the Expected value solution the average routing extra cost per client is of \$34 dollars whereas for the conservative solution \$22 dollars. This confirms that the conservative solution provides not only a cheaper expected routing cost value, but also the costs of serving an extra client are cheaper. 
 
\subsubsection{Total daily operation costs, and performance}\label{sss:dailyExpecctedCostst}

Table \ref{tab:metrics} we present several metrics that need to be considered when choosing a solution. Row 1 reports the total daily expected costs for both solutions (in dollars); row 2 reports the average vehicle utilization; row 3 reports the average daily CO2 emissions (in tons/day); and row 4 presents the expected service level for both solutions.

\begin{table}[h]
\centering\small
\caption{Expected value solution vs Conservative solution}
\label{tab:metrics}
\begin{tabular}{lrr}
\hline
                             & Expected value & Conservative \\ \hline
Daily costs                  & 3,665.5         & 3,733         \\
Vehicle utilization          & 92.20\%        & 75.60\%      \\
CO2 emissions  & 3.60            & 2.80          \\
Service level                & 99.70\%        & 100\%        \\ \hline
\end{tabular}
\end{table}

\begin{figure}
    \centering
    \subfloat[Vehicle utilization for the Expected value solution\label{fig:vehUtEV}]{\includegraphics[width=0.5\textwidth]{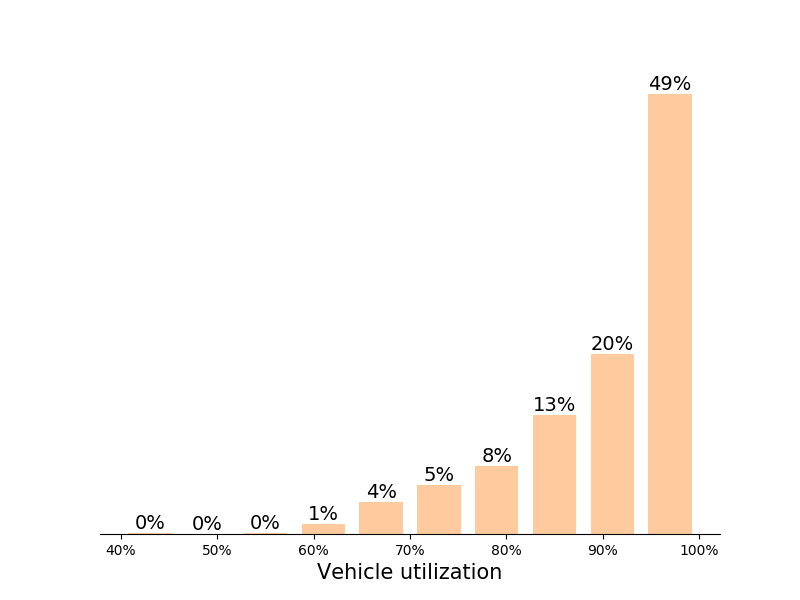}}
    \subfloat[Vehicle utilization for the Conservative solution\label{fig:vehUtC}]{\includegraphics[width=0.5\textwidth]{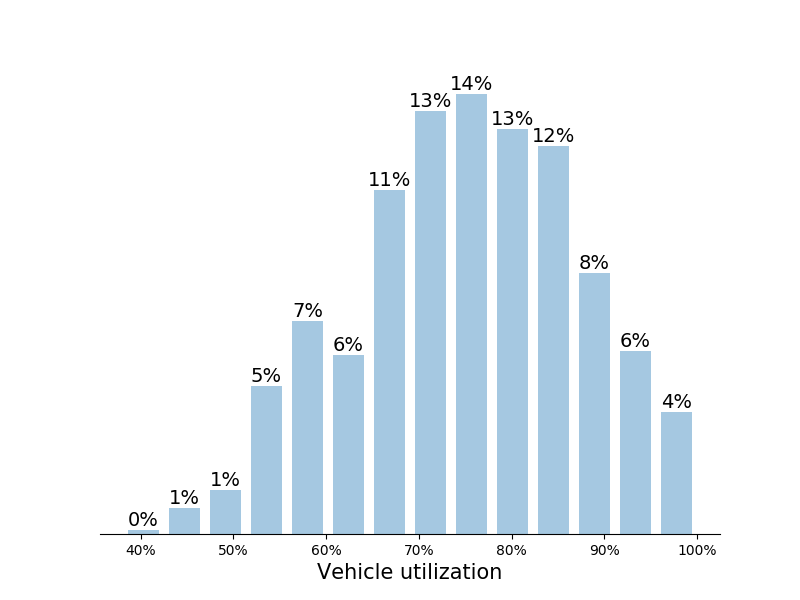}}
    \\
    \caption{Vehicle utilization}
    \label{fig:metrics}
\end{figure}

With the description of the routing costs obtained in \S \ref{sss:dailyExpecctedCostst}, and considering the daily maintenance costs, we can compare the total daily expected costs. For this task we performed a t--test and we found with a confidence of 95\%, and a p--value of 0.04 that there is enough statistical evidence that shows that the expected value on the total daily operation costs of the conservative solution is greater than with the expected value solution by \$67.5 dollars. Additionally, as shown in Figure \ref{fig:vehUtEV}, the average vehicle utilization of the expected value solution is 92.2\%, and in almost 50\% of the instances this solution utilizes all of its fleet (100\% vehicle utilization). On the other hand, the conservative solution has an average vehicle utilization of 75.5\%, and only 4\% of the instances utilize all the 21 vehicles placed (see Figure \ref{fig:vehUtC}). This meaning, on one hand, that the expected value solution uses the resources in a more efficient way, but also, that with this solution there is no slack for responding to unexpected high demand scenarios. Another important metric to consider is the expected daily CO2 emissions that for the conservative solution are 0.8 tons per day less than with the expected value solution, so in environmental terms, the conservative solution could be more appropriate. Finally it is really important to consider that expected value solution does not fulfill the total demand in 54 of the 730 days tested, whereas the conservative solution is able to serve all the clients in all the tested days.

\subsection{Computational times}\label{ss:compTimes}

Before the automation of the decision of the number of depots for each of the alternatives and the implementation of many acceleration strategies, we ran the algorithm many times with a fixed parameter for the number of depots to place (e.g. 1 to 9 depots). Table 3 \ref{tab:comptimes} compares the computational times of the algorithm for the expected value and conservative solutions using both the automated and non-automated approaches broken down by phases: warm start, training, facility location, set covering and route generators. Column 1 and 2 report times in seconds for the expected value solution in the non-automated and automated approaches. Column 3 reports the speedup of the automated approach with respect of the non-automated approach for the expected value solution. Columns 4 and 5 report times in seconds for the conservative solution in the non-automated and automated approaches. Finally, column 6 report the speedup of the automated approach with respect of the non-automated approach for the conservative approach.  

\begin{table}[htbp]
  \centering
  \setlength{\tabcolsep}{3pt}
  \small
  \caption{Computational times for the expected value and conservative instances (in seconds)}
    \begin{tabular}{lrrrrrr}
    \hline
          & \multicolumn{3}{c}{Expected value} & \multicolumn{3}{c}{Conservative}\\
\cline{2-7}          & \multicolumn{1}{c}{Non-automated} & \multicolumn{1}{c}{Automated} & \multicolumn{1}{c}{Speedup} & \multicolumn{1}{c}{Non-automated} & \multicolumn{1}{c}{Automated} & \multicolumn{1}{c}{Speedup}\\
    \hline
    Warm start &                        0.56  &                      0.09  &             6.01  &                      1.05  &                      0.22  &             4.79 \\
    Training &                1,879.19  &              1,879.19  &             1.00  &              4,940.51  &              4,940.51  &             1.00  \\
    Facility location &                    738.13  &              2,363.49  &             0.31  &              1,169.33  &              3,675.89  &             0.32  \\
    Set covering &                1,315.76  &                 655.14  &             2.01  &              2,663.36  &              1,102.76  &             2.42  \\
    Route generators &              11,782.63  &              1,153.61  &          10.21  &            49,891.27  &            11,147.37  &             4.48  \\
    \hline 
    Total Time &              15,716.27  &              6,051.53  &             2.60  &            58,665.51  &            20,866.74  &             2.81 \\
    \hline
    \end{tabular}%
  \label{tab:comptimes}%
\end{table}%

The automated algorithm is approximately 2.7 times faster than the non-automated approach. The phase with the greatest speedup is the warm start and after it the route generators phase. The later is 4.48 times faster in the conservative solution and 10.21 times faster in the expected value solution than the non-automated approach. The total execution time for the end user is about 5.8 hours (20,866.74 s) in the conservative solution and 1.7 hours (15,716 s) in the expected value solution. Meaning that in less than 6 hours the decision makers can count with a solution for the strategic decision of depot location and fleet allocation. 

\begin{figure}[h]
    \centering
    \includegraphics[width=0.7\textwidth]{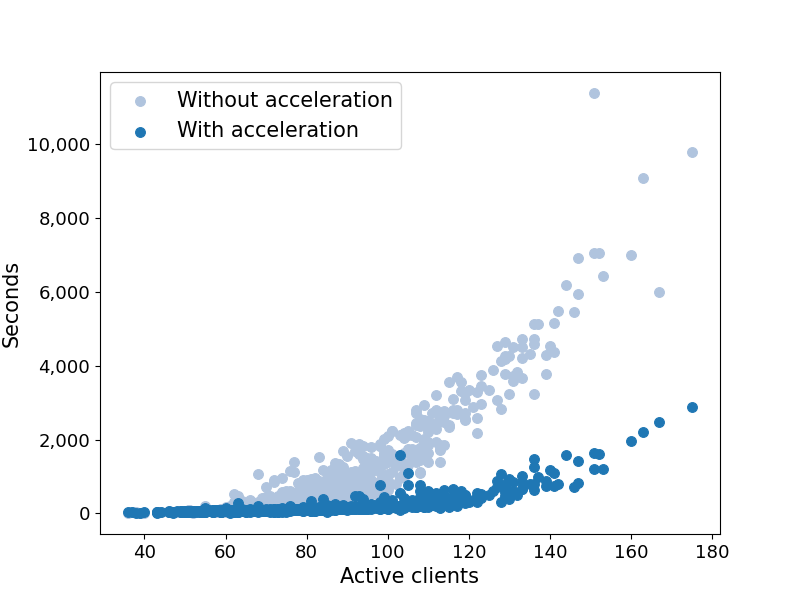}
    \caption{Computational times for solving MDVRP's with and without our acceleration strategies (in seconds).}
    \label{fig:MDVRPCompTimes}
\end{figure}

Finally, Figure \ref{fig:MDVRPCompTimes} summarizes the computational times of the 730 solved instances solved with acceleration strategies (in dark blue) and without acceleration strategies (in light blue). In the horizontal axis is the number of daily active clients and in the vertical axis is the computational time. Note that the relation between the scale of the problem and the time to solve grows exponentially when solving without acceleration strategies. When acceleration techniques are not used, the proposed methodology solves the 65\% of the instances in less than 30 minutes, and 30\% of the instances in less than 10 minutes. With acceleration strategies, the proposed methodology solves 75\% of the instances in less than 3 minutes. The most difficult instance which had 175 active clients and a total demand of 1741 units took 40 minutes to solve with acceleration strategies, meaning that the model is appropriate to suggest a routing plan to support operational decisions.

\section{Conclusions}\label{s:Conclusions}

The methodology consists of a two-stage stochastic program with a facility location problem in the first stage and a routing generation scheme in the second stage. The route generation is performed via a column generation approach, where the master problem select routes via a set covering model and the auxiliary problem receives pricing information from the master problem to build routes in a transformed graph using a ESPPRC-R formulation. The routing problems give information back to the facility location in the first stage in an iterative fashion by means of cuts combine classic continuous bender cuts with incentives to binary variables. We consider that our main contribution is the articulation of all these elements into one complete solution scheme. 

The generation of scenarios for the second stage problems is done by means of simulation with the parameters of the adjusted in the descriptive models. To select the representative scenarios, the risk profile of the decision makers is considered and two alternatives are proposed: expected value and conservative case. In the conservative case the goal is to increase service level using robust solutions and providing worst-case. However, incorporating too many scenarios could compromise obtaining a practical solution in a reasonable time. Therefore, the simulation based approach to obtain scenarios serves well if the statistical properties of the dataset are explored beforehand. 

Using a two-stage stochastic model scheme is advantageous because the first stage receives some recourse information from the second stage scenarios and takes efficient decisions on how to install resources as depots and vehicle fleet strategically. The efficient configuration of depots and vehicles minimizes the costs of the daily operation, offering a service level of 100\% if the conservative solution is selected. We provide in the report relevant metrics of the solutions for each alternative so decision makers can select the solution to use. 

The results were framed in a decision support system built in AIMMS. The tool receives a demand file with the observations to solve for the locations in the state of Pennsylvania and provides three options to solve for those observations: (i) petal recycled heuristic described in \S \ref{sss:AcelSP} and (ii) column generation optimally with the pulse algorithm described in \S \ref{ss:SecondStage} (iii) just one iteration of the column generation. The second and the last option run the same model with the difference that (i) runs until the column generation does not find more improvements, which will take more time and (ii) runs the auxiliary problem (route generator) and the master problem (route selection) of the column generation one time, so the user can get solutions in less than 40 minutes and see how the algorithm improves in each iteration. For each run, the decision support tool shows the depots and the clients over the map, the number of vehicles used, the cost of routing that particular observation. If a depot is selected then the routes assigned to that depot appear. The interface and results display were designed to provide plans for the daily operation.

Future work can focus on reducing computational times. Experiments with heuristic approaches can provide insightful information about the trade off between time reduction and quality loss in the solutions.

\bibliographystyle{spbasic}      
\bibliography{VFSPRP.bib}

\appendix
\section{Algorithms}\label{Ap:Algorithms}
\subsection{Column Generation algorithm}

Algorithm \ref{alg:ColumnGeneration} outlines the specific CG followed to solve the MDVRP. Given a specific first-stage solution $\hat{\bm{x}},\hat{\bm{y}}$ and second-stage problem data $\bm{\xi}(\omega)$ the procedure estimates the second-stage value $\mathcal{Q}(\hat{\bm{x}},\hat{\bm{y}},\bm{\xi}(\omega))$. Line 1 sets the stopping condition to \texttt{False}; line 2 starts repeating lines 3--11 until no route with negative reduced cost is found by the route generator (see \S \ref{ss:RouteGenerator}); line 4 solves the RRMP and recovers the dual vectors $\bm{\pi}$ and $\bm{\lambda}$ of constraints \ref{eq:SPMset_part} and \ref{eq:SPMcap}, respectively; lines 5--11 generate routes starting and ending at each $k\in\hat{\Hh}$ (line 6); if the generated route has negative reduced cost, the route is added to the set $\Rout_k(\omega)$ (line 7 and 8); finally line 13 solves the RMP with the previously generated set of routes $\Rout(\omega)$ and returns the second stage value for this first-stage solution $\hat{\bm{x}},\hat{\bm{y}}$ and second-stage problem data $\bm{\xi}(\omega)$, $\mathcal{Q}(\hat{\bm{x}},\hat{\bm{y}},\bm{\xi}(\omega))$.

\begin{algorithm}[hbt]
  \begin{algorithmic}[1]
    \Require {First stage solution $\hat{\bm{x}},\hat{\bm{y}}$ and second-stage problem data $\bm{\xi}(\omega)$}
    \Ensure {Second-stage value $\mathcal{Q}(\hat{\bm{x}},\hat{\bm{y}},\bm{\xi}(\omega))$}
    \State{Stop$\leftarrow\texttt{False}$}
    \Repeat
        \State{Stop$\leftarrow\texttt{True}$}
        \State{$\bm{\lambda,\pi}\leftarrow$Solve RRMP}
        \For{$h\in\hat{\Hh}$}
            \State{$\bm{\lambda,\pi}\leftarrow$Solve RRMP}
            \State{$r\leftarrow$\texttt{RouteGenerator}($h,\lambda_h,\bm{\pi}$)}\Comment{see \S \ref{ss:RouteGenerator}}
            \If{$\bar{c}_r<0$}
                \State{$\Rout_k(\omega)\leftarrow \Rout_k(\omega)\cup \{r\}$}
                \State{Stop$\leftarrow\texttt{False}$}
            \EndIf
        \EndFor
    \Until{Stop}
    \State{$\mathcal{Q}(\hat{\bm{x}},\hat{\bm{y}},\bm{\xi}(\omega))\leftarrow$Solve RMP}
    
    \Return {$\mathcal{Q}(\hat{\bm{x}},\hat{\bm{y}},\bm{\xi}(\omega))$}
  \end{algorithmic}
  \caption{Column Generation for the MDVRP.}
  \label{alg:ColumnGeneration}
\end{algorithm}

\subsection{L-shaped algorithm}

Line 1 sets the number of iterations $t$ to zero, the upper bound $Ub$ to infinity and the lower bound $Lb$ to minus infinity; line 2 starts repeating the lines 3--12 until the optimality gap, namely $\frac{Ub-Lb}{Ub}$, is below our optimality tolerance $\varepsilon$; line 3 solves the MP and recovers the first stage solutions $\hat{\bm{x}}$ and $\hat{\bm{y}}$, the estimations of the second-stage value $\bm{\eta}$, and the current estimation of the expected total cots $\nu^{t^*}$; line 4 sets the lower bound to the maximum between the current lower bound $Lb$ and the current estimation of the expected total cots $\nu^{t^*}$; lines 5--10 solve for each $s\in\Scen$ a MDVRP with the first stage solutions $\hat{\bm{x}}$ and $\hat{\bm{y}}$ and its respective second-stage problem data $\bm{\xi}(s)$ (line 6), if the MP underestimates $Q(\hat{\bm{x}},\hat{\bm{y}},\bm{\xi}(s))$ (i.e., $\eta_s<\eta_s^{t^*}$), a cut \eqref{eq:BMPOpt} is added to the MP; line 11 updates the upper bound by setting it to the minimum between the current upper bound $Ub$ and the sum of the current first-stage costs and the expected value of the second-stage for the first stage solution $\hat{\bm{x}}$ and $\hat{\bm{y}}$, namely $\sum_{s\in\Scen}p_s\eta_s^{t^*}$; line 12 updates the number the iteration; and finally an optimal solution $(\bm{x^*,y^*,\eta^*})$ is returned.
\begin{algorithm}[hbt]
  \begin{algorithmic}[1]
    \Require {Optimality gap $\varepsilon$}
    \Ensure {$\varepsilon$-Optimal solution for the VFSPRP-SC $\bm{x}^*,\bm{y}^*,\bm{\eta}^*$}
    \State{$t,Ub,Lb\leftarrow 0,\infty,-\infty$}
    \While { $Ub-Lb>\varepsilon \cdot Ub$}
        \State{$\hat{\bm{x}},\hat{\bm{y}},\bm{\eta},\nu^{t^*} \leftarrow$ Solve MP}
        \State{$Lb\leftarrow\max\{Lb,\nu^{t^*}\}$}
        \For{$s\in\Scen$}
            \State{$\eta^{t^*}_s\leftarrow \texttt{ColumGeneration}(\hat{\bm{x}},\hat{\bm{y}},\bm{\xi}(s))$}\Comment{Algorithm \ref{alg:ColumnGeneration}}
            \If{$\eta_s<\eta^{t^*}_s$}
                \State{Add cut (\ref{eq:BMPOpt}) to the MP}
            \EndIf
        \EndFor
        \State{$Ub\leftarrow\min\left\{Ub,\nu^{t^*}-\sum_{s\in\Scen}p_s(\eta_s-\eta^{t^*}_s)\right\}$}
        \State{$t\leftarrow t+1$}
    \EndWhile
    
    \Return {$\bm{x}^*,\bm{y}^*,\bm{\eta}^*$}
  \end{algorithmic}
  \caption{Multicut L-Shaped Algorithm for the VFSPRP-SC.}
  \label{alg:LShaped}
\end{algorithm}	

\begin{algorithm}[hbt]
  \begin{algorithmic}[1]
    \Require {Optimality gap $\varepsilon$}
    \Ensure {$\varepsilon$-Optimal solution for the VFSPRP-SC $\bm{x}^*,\bm{y}^*,\bm{\eta}^*$}
    
    \State{$\hat{\bm{x}},\hat{\bm{y}},\bm{\eta},\nu^{t^*} \leftarrow$ \texttt{warmStart}()}
    
    \State{$t,Ub(0),Lb(0)\leftarrow 0,\infty,-\infty$,\texttt{ True}}
    \State{$\text{RouteGen}\leftarrow\texttt{True}$}
    
    \While { $Ub(t)-Lb(t)>\varepsilon \cdot Ub(t)$}
        \State{$\hat{\bm{x}},\hat{\bm{y}},\bm{\eta},\nu^{t^*} \leftarrow$ Solve MP}
        \State{$Lb(t+1)\leftarrow\max\{Lb(t),\nu^{t^*}\}$}
        \For{$s\in\Scen$}
        \If{RouteGen}
            \State{$\eta^{t^*}_s\leftarrow \texttt{ColumGeneration}(\hat{\bm{x}},\hat{\bm{y}},\bm{\xi}(s))$}
        \Else
            \State{$\eta^{t^*}_s\leftarrow \texttt{SolveRRMP}(\hat{\bm{x}},\hat{\bm{y}},\bm{\xi}(s))$}\Comment{Algorithm \ref{alg:ColumnGeneration}}
        \EndIf
            
            \If{$\eta_s<\eta^{t^*}_s$}
                \State{Add cut (\ref{eq:BMPOpt}) to the MP}
            \EndIf
        \EndFor
        \State{$Ub(t+1)\leftarrow\min\left\{Ub(t),\nu^{t^*}-\sum_{s\in\Scen}p_s(\eta_s-\eta^{t^*}_s)\right\}$}
        \If{$Ub(t+1) < Ub(t)$}
            \State{RouteGen$\leftarrow$\texttt{True}}
        \Else
            \State{RouteGen$\leftarrow$\texttt{False}}
        \EndIf
        \State{$t\leftarrow t+1$}
    \EndWhile
    
    \Return {$\bm{x}^*,\bm{y}^*,\bm{\eta}^*$}
  \end{algorithmic}
  \caption{Multicut L-Shaped Algorithm for the VFSPRP-SC.}
  \label{alg:LShaped1}
\end{algorithm}

\subsection{Petal recycler algorithm}
Algorithm \ref{alg:PetalRecycler} formally outlines the procedure followed to map a route $r\in\Rout^s_h$ into a route $r'\in\Rout^{s'}_h$. Line 1 initializes the new route $r'\in\Rout^{s'}_h$, its load $Q_{r'}$ and its current time consumption $T_{r'}$; line 2 saves the $\ell-1$-th node in $r'$; line 3 starts replacing every node $v_i \in r$ for nodes in $\D(s)$ by taking as candidate the node $v_\ell$ with minimum $\bar{c}(v_{\ell-1},v_i,v_\ell)$ (line 4), if the load constraint is met (line 5) and the time windows constraints are met (line 6), the node $v_\ell$ is added to route $r'$ (line 7) and the $\ell-1$-th node in $r'$ is updated (line 8); finally it returns a high quality, feasible route $r'\in\Rout^{s'}_h$.

\begin{algorithm}[hbt]
  \begin{algorithmic}[1]
    \Require {Route for specific second stage problem data $\bm{\xi}(s)$, $r\in\Rout^s_h$; second stage problem data for which $r$ wants to be mapped, $\bm{\xi}(s')$.}
    \Ensure {Feasible route for the second stage problem data $\bm{\xi}^{s'}$, $r'\in\Rout^{s'}_h$.}
    \State{$r',Q_{r'},T_{r'}\leftarrow\{h\},0,0$}
    \State{$v_{\ell-1}\leftarrow h$}
    \For{$v_i \in r $}
        \State{$v_\ell \leftarrow \argmin\limits_{v_k\in\D(s)}\{\bar{c}(v_{\ell-1},v_i,v_k)\}$}
        \If{$Q_{r'}+d^s_\ell<Q$}
            \If{$w_1(\ell)\le T_{r'}+\tau^{sp}(v_{\ell-1},v_{\ell})\le w_2(\ell)$}
                \State{$r'\leftarrow r'\cup \{v_\ell\}$}
                \State{$v_{\ell-1}\leftarrow v_\ell$}
            \EndIf
        \EndIf
    \EndFor
    
    \Return {$r'$}
  \end{algorithmic}
  \caption{Petal Recycler.}
  \label{alg:PetalRecycler}
\end{algorithm}

\section{Models}\label{Ap:models}

\subsection{Scenario selection model}\label{Aps:scensel}

Let $\mathcal{T}$ be the set of realizations classified to be part of a scenario. Let $\mathcal{K}$ be the set of districts and $\N_k$ the subset of clients (nodes) in each district $k \in \mathcal{K}$. $d^{t}_{\N_k}$ be the parameter of the aggregated demand of $\N_k$ in realization $t\in \mathcal{T}$ of $\tilde{d}_{\N'} \triangleq \sum_{i \in \N'}\tilde{d}_{i}$ where $\tilde{d}_{i}$ is the demand of client $i \in \N'$. The model decides with variable $x^{t}$ whether a realization $t \in \mathcal{T}$ is representative of the scenario. With auxiliary variable $z$, the model captures the minimum demand between districts.

\begin{maxi!}|s|
    {}{z \protect\label{eq:ScenGObj}}
    {\label{prob:ScenG}}{}
    \addConstraint{\sum_{t\in \mathcal{T}}d^{t}_{\N_k}x^{t}}{\geq z}{\quad \forall k \in \mathcal{K} \protect\label{eq:ScenGMaxMinDav}}
    \addConstraint{\sum_{i\in \mathcal{T}}x^{t}}{=1}{\ \protect\label{eq:ScenGOneR}}
    \addConstraint{x^{t}}{\in \{0,1\},}{\quad \forall t\in \mathcal{T} \protect\label{eq:ScenGNatx}}
    \addConstraint{z}{\in \mathbb{Z}_{+}^{1}.}{\protect\label{eq:ScenGNatz}}
\end{maxi!}\\

For each scenario $s\in \Scen$ model \ref{prob:ScenG} should be solved. Objective function \ref{eq:ScenGObj} maximize auxiliary variable $z$ which captures with \ref{eq:ScenGMaxMinDav} the minimum demand within the districts in the selected realization. The group of constraints \ref{eq:ScenGOneR} guarantees that only one realization is selected from the sample. Finally, constraints \ref{eq:ScenGNatx} and \ref{eq:ScenGNatz} define the nature of variables.
\end{document}